\documentclass{gtart}

\def\ifplaintex{\expandafter\ifx\csname documentclass\endcsname\relax}

\def\gtp{{\mathsurround=0pt\it $\cal G\mskip-2mu$eometry \&\ 
$\cal T\!\!$opology $\cal P\!$ublications}}  

\def\recd{{\small Received:\qua\receiveddate\ifx\reviseddate\relax
\else\qquad Revised:\qua\reviseddate\fi\par}} 

\def\lognumber#1{\def\thelognumber{#1}}
\def\volumenumber#1{\def\thevolumenumber{#1}}
\def\volumeyear#1{\def\thevolumeyear{#1}}
\def\papernumber#1{\def\thepapernumber{#1}}
\def\pagenumbers#1#2{\def\startpage{#1}\def\finishpage{#2}}
\def\published#1{\def\publishdate{#1}}

\def\received#1{\def\receiveddate{#1}}

\def\accepted#1{\def\accepteddate{#1}}
\def\asciititle#1{\def\theasciititle{#1}}
\def\covertitle#1{\def\thecovertitle{#1}}

\def\asciiaddress#1{\def\theasciiaddress{#1}}

\long\def\asciiabstract#1{\long\def\theasciiabstract{#1}}
\def\asciikeywords#1{\def\theasciikeywords{#1}}

\let\\\par\let\thelognumber\relax\let\thevolumenumber\relax
\let\thepapernumber\relax\let\thevolumeyear\relax\let\startpage\relax
\let\finishpage\relax\let\publishdate\relax\let\receiveddate\relax
\let\reviseddate\relax\let\accepteddate\relax\let\theasciititle\relax
\let\thecovertitle\relax\let\theasciiauthors\relax\let\theasciiaddress\relax
\let\theasciiabstract\relax\let\theasciikeywords\relax

\let\theasciiemail\relax

\ifplaintex
\font\logobig=cmssbx10 scaled 3836
\font\logomed=cmssbx10 scaled 2557
\else
\font\logobig=cmssbx10 scaled 4200
\font\logomed=cmssbx10 scaled 2800
\fi

\long\def\makeagttitle{   
\count0=\startpage
\agt\hfill      
\hbox to 45truept{\vbox to 0pt{\vglue -13truept{\logomed A\kern -.37em{\logobig 
T}\kern -.38em G}\vss}\hss}
\break
{\small Volume \thevolumenumber\ (\thevolumeyear)
\startpage--\finishpage\nl
Published: \publishdate}

\vglue .25truein

{\parskip=0pt\leftskip 0pt plus
1fil\def\\{\par\smallskip}{\Large\bf\thetitle}\par\medskip} \vglue
0.05truein

{\parskip=0pt\leftskip 0pt plus 1fil\def\\{\par}{\sc\theauthors}
\par\medskip}%
 
\vglue 0.03truein 

{\small\leftskip 25truept\rightskip 25truept{\bf Abstract}\stdspace\theabstract

{\bf AMS Classification}\stdspace\theprimaryclass
\ifx\thesecondaryclass\relax\else; \thesecondaryclass\fi\par
{\bf Keywords}\stdspace \thekeywords\par}\vglue 7truept

}   

\ifplaintex
\font\phead=cmsl9 scaled 950
\font\pnum=cmbx10 scaled 913
\font\pfoot=cmsl9 scaled 950
\headline{\vbox to 0pt{\vskip -4.5mm\line{\small\phead\ifnum
\count0=\startpage ISSN 1472-2739 (on-line) 1472-2747 (printed)
\hfill {\pnum\folio}\else\ifodd\count0\def\\{ }%
\ifx\theshorttitle\relax\thetitle\else\theshorttitle\fi\hfill{\pnum\folio}
\else\def\\{ and }{\pnum\folio}\hfill\ifx\theshortauthors\relax\theauthors
\else\theshortauthors\fi\fi\fi}\vss}}
\footline{\vbox to 0pt{\vglue 0mm\line{\small\pfoot\ifnum\count0=\startpage
\copyright\ \gtp\hfill\else
\agt, Volume \thevolumenumber\ (\thevolumeyear)\hfill\fi}\vss}}
\else
\font=cmsl9 scaled 1050
\font\lnum=cmbx10 
\font=cmsl9 scaled 1050
\makeatletter
\def\@oddhead{{\small\lhead\ifnum\count0=\startpage ISSN 1472-2739 
(on-line) 1472-2747 (printed)\hfill {\lnum\number\count0}\else\ifodd\count0
\def\\{ }\ifx\theshorttitle\relax \thetitle \else\theshorttitle\fi\hfill
{\lnum\number\count0}\else\def\\{ and }{\lnum\number\count0}
\hfill\ifx\theshortauthors\relax 
\theauthors\else\theshortauthors\fi\fi\fi}}\def\@evenhead{\@oddhead}
\def\@oddfoot{\small\lfoot\ifnum\count0=\startpage\copyright\ \gtp\hfill\else
\agt, Volume \thevolumenumber\ (\thevolumeyear)\hfill\fi}
\def\@evenfoot{\@oddfoot}
\makeatother
\fi
\let\maketitlepage\makeagttitle
\let\makeshorttitle\maketitlepage
\let\maketitle\maketitlepage

\newwrite\gtoutfile
\long\gdef\makeheadfile{  
{\def\\{, }\def\s{ }
\immediate\openout\gtoutfile head.xxx
\immediate\write\gtoutfile{To: math@arxiv.org}
\immediate\write\gtoutfile{Subject: put OR rep NNNNN:ppppp}
\immediate\write\gtoutfile{--text follows this line--}
\immediate\write\gtoutfile{Proxy-for: \ifx\theasciiauthors\relax
\theauthors\else\theasciiauthors\fi\s<\ifx\theasciiemail\relax\theemail\else\theasciiemail\fi>}
\immediate\write\gtoutfile{\noexpand\\}
\immediate\write\gtoutfile{Authors: \ifx\theasciiauthors\relax
\theauthors\else\theasciiauthors\fi}
{\def\\{ }\immediate\write\gtoutfile{Title: \ifx\theasciititle\relax
\thetitle\else\theasciititle\fi}}
\immediate\write\gtoutfile{Subj-class: GT or SG, GR etc}
\immediate\write\gtoutfile{MSC-class: \theprimaryclass\ifx\thesecondaryclass\relax\else, \thesecondaryclass\fi}
\immediate\write\gtoutfile{Journal-ref: Algebr. Geom. Topol. \thevolumenumber\s
(\thevolumeyear) \startpage-\finishpage}
\immediate\write\gtoutfile{Comments: Published by Algebraic and
Geometric Topology at}
\immediate\write\gtoutfile{\s\s\s  http://www.maths.warwick.ac.uk/agt/AGTVol\thevolumenumber/agt-\thevolumenumber-\thepapernumber.abs.html}
\immediate\write\gtoutfile{\noexpand\\}
\immediate\write\gtoutfile{}
\ifx\theasciiabstract\relax
\immediate\write\gtoutfile{\theabstract}\else
\immediate\write\gtoutfile{\theasciiabstract}\fi
\immediate\write\gtoutfile{}
\immediate\write\gtoutfile{\noexpand\\}
\immediate\write\gtoutfile{}
\immediate\closeout\gtoutfile}}  

\def\maketitlepage{\makeagttitle\makeheadfile}
\let\makeshorttitle\maketitlepage
\let\maketitle\maketitlepage

\def\ifplaintex{\expandafter\ifx\csname documentclass\endcsname\relax}

\def\gtp{{\mathsurround=0pt\it $\cal G\mskip-2mu$eometry \&\ 
$\cal T\!\!$opology $\cal P\!$ublications}}  

\def\recd{{\small Received:\qua\receiveddate\ifx\reviseddate\relax
\else\qquad Revised:\qua\reviseddate\fi\par}} 

\def\lognumber#1{\def\thelognumber{#1}}
\def\volumenumber#1{\def\thevolumenumber{#1}}
\def\volumeyear#1{\def\thevolumeyear{#1}}
\def\papernumber#1{\def\thepapernumber{#1}}
\def\pagenumbers#1#2{\def\startpage{#1}\def\finishpage{#2}}
\def\published#1{\def\publishdate{#1}}

\def\received#1{\def\receiveddate{#1}}

\def\accepted#1{\def\accepteddate{#1}}
\def\asciititle#1{\def\theasciititle{#1}}
\def\covertitle#1{\def\thecovertitle{#1}}

\def\asciiaddress#1{\def\theasciiaddress{#1}}

\long\def\asciiabstract#1{\long\def\theasciiabstract{#1}}
\def\asciikeywords#1{\def\theasciikeywords{#1}}

\let\\\par\let\thelognumber\relax\let\thevolumenumber\relax
\let\thepapernumber\relax\let\thevolumeyear\relax\let\startpage\relax
\let\finishpage\relax\let\publishdate\relax\let\receiveddate\relax
\let\reviseddate\relax\let\accepteddate\relax\let\theasciititle\relax
\let\thecovertitle\relax\let\theasciiauthors\relax\let\theasciiaddress\relax
\let\theasciiabstract\relax\let\theasciikeywords\relax

\let\theasciiemail\relax

\ifplaintex
\font\logobig=cmssbx10 scaled 3836
\font\logomed=cmssbx10 scaled 2557
\else
\font\logobig=cmssbx10 scaled 4200
\font\logomed=cmssbx10 scaled 2800
\fi

\long\def\makeagttitle{   
\count0=\startpage
\agt\hfill      
\hbox to 45truept{\vbox to 0pt{\vglue -13truept{\logomed A\kern -.37em{\logobig 
T}\kern -.38em G}\vss}\hss}
\break
{\small Volume \thevolumenumber\ (\thevolumeyear)
\startpage--\finishpage\nl
Published: \publishdate}

\vglue .25truein

{\parskip=0pt\leftskip 0pt plus
1fil\def\\{\par\smallskip}{\Large\bf\thetitle}\par\medskip} \vglue
0.05truein

{\parskip=0pt\leftskip 0pt plus 1fil\def\\{\par}{\sc\theauthors}
\par\medskip}%
 
\vglue 0.03truein 

{\small\leftskip 25truept\rightskip 25truept{\bf Abstract}\stdspace\theabstract

{\bf AMS Classification}\stdspace\theprimaryclass
\ifx\thesecondaryclass\relax\else; \thesecondaryclass\fi\par
{\bf Keywords}\stdspace \thekeywords\par}\vglue 7truept

}   

\ifplaintex
\font\phead=cmsl9 scaled 950
\font\pnum=cmbx10 scaled 913
\font\pfoot=cmsl9 scaled 950
\headline{\vbox to 0pt{\vskip -4.5mm\line{\small\phead\ifnum
\count0=\startpage ISSN 1472-2739 (on-line) 1472-2747 (printed)
\hfill {\pnum\folio}\else\ifodd\count0\def\\{ }%
\ifx\theshorttitle\relax\thetitle\else\theshorttitle\fi\hfill{\pnum\folio}
\else\def\\{ and }{\pnum\folio}\hfill\ifx\theshortauthors\relax\theauthors
\else\theshortauthors\fi\fi\fi}\vss}}
\footline{\vbox to 0pt{\vglue 0mm\line{\small\pfoot\ifnum\count0=\startpage
\copyright\ \gtp\hfill\else
\agt, Volume \thevolumenumber\ (\thevolumeyear)\hfill\fi}\vss}}
\else
\font=cmsl9 scaled 1050
\font\lnum=cmbx10 
\font=cmsl9 scaled 1050
\makeatletter
\def\@oddhead{{\small\lhead\ifnum\count0=\startpage ISSN 1472-2739 
(on-line) 1472-2747 (printed)\hfill {\lnum\number\count0}\else\ifodd\count0
\def\\{ }\ifx\theshorttitle\relax \thetitle \else\theshorttitle\fi\hfill
{\lnum\number\count0}\else\def\\{ and }{\lnum\number\count0}
\hfill\ifx\theshortauthors\relax 
\theauthors\else\theshortauthors\fi\fi\fi}}\def\@evenhead{\@oddhead}
\def\@oddfoot{\small\lfoot\ifnum\count0=\startpage\copyright\ \gtp\hfill\else
\agt, Volume \thevolumenumber\ (\thevolumeyear)\hfill\fi}
\def\@evenfoot{\@oddfoot}
\makeatother
\fi
\let\maketitlepage\makeagttitle
\let\makeshorttitle\maketitlepage
\let\maketitle\maketitlepage

\newwrite\gtoutfile
\long\gdef\makeheadfile{  
{\def\\{, }\def\s{ }
\immediate\openout\gtoutfile head.xxx
\immediate\write\gtoutfile{To: math@arxiv.org}
\immediate\write\gtoutfile{Subject: put OR rep NNNNN:ppppp}
\immediate\write\gtoutfile{--text follows this line--}
\immediate\write\gtoutfile{Proxy-for: \ifx\theasciiauthors\relax
\theauthors\else\theasciiauthors\fi\s<\ifx\theasciiemail\relax\theemail\else\theasciiemail\fi>}
\immediate\write\gtoutfile{\noexpand\\}
\immediate\write\gtoutfile{Authors: \ifx\theasciiauthors\relax
\theauthors\else\theasciiauthors\fi}
{\def\\{ }\immediate\write\gtoutfile{Title: \ifx\theasciititle\relax
\thetitle\else\theasciititle\fi}}
\immediate\write\gtoutfile{Subj-class: GT or SG, GR etc}
\immediate\write\gtoutfile{MSC-class: \theprimaryclass\ifx\thesecondaryclass\relax\else, \thesecondaryclass\fi}
\immediate\write\gtoutfile{Journal-ref: Algebr. Geom. Topol. \thevolumenumber\s
(\thevolumeyear) \startpage-\finishpage}
\immediate\write\gtoutfile{Comments: Published by Algebraic and
Geometric Topology at}
\immediate\write\gtoutfile{\s\s\s  http://www.maths.warwick.ac.uk/agt/AGTVol\thevolumenumber/agt-\thevolumenumber-\thepapernumber.abs.html}
\immediate\write\gtoutfile{\noexpand\\}
\immediate\write\gtoutfile{}
\ifx\theasciiabstract\relax
\immediate\write\gtoutfile{\theabstract}\else
\immediate\write\gtoutfile{\theasciiabstract}\fi
\immediate\write\gtoutfile{}
\immediate\write\gtoutfile{\noexpand\\}
\immediate\write\gtoutfile{}
\immediate\closeout\gtoutfile}}  

\def\maketitlepage{\makeagttitle\makeheadfile}
\let\makeshorttitle\maketitlepage
\let\maketitle\maketitlepage

\def\ifplaintex{\expandafter\ifx\csname documentclass\endcsname\relax}

\def\gtp{{\mathsurround=0pt\it $\cal G\mskip-2mu$eometry \&\ 
$\cal T\!\!$opology $\cal P\!$ublications}}  

\def\recd{{\small Received:\qua\receiveddate\ifx\reviseddate\relax
\else\qquad Revised:\qua\reviseddate\fi\par}} 

\def\lognumber#1{\def\thelognumber{#1}}
\def\volumenumber#1{\def\thevolumenumber{#1}}
\def\volumeyear#1{\def\thevolumeyear{#1}}
\def\papernumber#1{\def\thepapernumber{#1}}
\def\pagenumbers#1#2{\def\startpage{#1}\def\finishpage{#2}}
\def\published#1{\def\publishdate{#1}}

\def\received#1{\def\receiveddate{#1}}

\def\accepted#1{\def\accepteddate{#1}}
\def\asciititle#1{\def\theasciititle{#1}}
\def\covertitle#1{\def\thecovertitle{#1}}

\def\asciiaddress#1{\def\theasciiaddress{#1}}

\long\def\asciiabstract#1{\long\def\theasciiabstract{#1}}
\def\asciikeywords#1{\def\theasciikeywords{#1}}

\let\\\par\let\thelognumber\relax\let\thevolumenumber\relax
\let\thepapernumber\relax\let\thevolumeyear\relax\let\startpage\relax
\let\finishpage\relax\let\publishdate\relax\let\receiveddate\relax
\let\reviseddate\relax\let\accepteddate\relax\let\theasciititle\relax
\let\thecovertitle\relax\let\theasciiauthors\relax\let\theasciiaddress\relax
\let\theasciiabstract\relax\let\theasciikeywords\relax

\let\theasciiemail\relax

\ifplaintex
\font\logobig=cmssbx10 scaled 3836
\font\logomed=cmssbx10 scaled 2557
\else
\font\logobig=cmssbx10 scaled 4200
\font\logomed=cmssbx10 scaled 2800
\fi

\long\def\makeagttitle{   
\count0=\startpage
\agt\hfill      
\hbox to 45truept{\vbox to 0pt{\vglue -13truept{\logomed A\kern -.37em{\logobig 
T}\kern -.38em G}\vss}\hss}
\break
{\small Volume \thevolumenumber\ (\thevolumeyear)
\startpage--\finishpage\nl
Published: \publishdate}

\vglue .25truein

{\parskip=0pt\leftskip 0pt plus
1fil\def\\{\par\smallskip}{\Large\bf\thetitle}\par\medskip} \vglue
0.05truein

{\parskip=0pt\leftskip 0pt plus 1fil\def\\{\par}{\sc\theauthors}
\par\medskip}%
 
\vglue 0.03truein 

{\small\leftskip 25truept\rightskip 25truept{\bf Abstract}\stdspace\theabstract

{\bf AMS Classification}\stdspace\theprimaryclass
\ifx\thesecondaryclass\relax\else; \thesecondaryclass\fi\par
{\bf Keywords}\stdspace \thekeywords\par}\vglue 7truept

}   

\ifplaintex
\font\phead=cmsl9 scaled 950
\font\pnum=cmbx10 scaled 913
\font\pfoot=cmsl9 scaled 950
\headline{\vbox to 0pt{\vskip -4.5mm\line{\small\phead\ifnum
\count0=\startpage ISSN 1472-2739 (on-line) 1472-2747 (printed)
\hfill {\pnum\folio}\else\ifodd\count0\def\\{ }%
\ifx\theshorttitle\relax\thetitle\else\theshorttitle\fi\hfill{\pnum\folio}
\else\def\\{ and }{\pnum\folio}\hfill\ifx\theshortauthors\relax\theauthors
\else\theshortauthors\fi\fi\fi}\vss}}
\footline{\vbox to 0pt{\vglue 0mm\line{\small\pfoot\ifnum\count0=\startpage
\copyright\ \gtp\hfill\else
\agt, Volume \thevolumenumber\ (\thevolumeyear)\hfill\fi}\vss}}
\else
\font=cmsl9 scaled 1050
\font\lnum=cmbx10 
\font=cmsl9 scaled 1050
\makeatletter
\def\@oddhead{{\small\lhead\ifnum\count0=\startpage ISSN 1472-2739 
(on-line) 1472-2747 (printed)\hfill {\lnum\number\count0}\else\ifodd\count0
\def\\{ }\ifx\theshorttitle\relax \thetitle \else\theshorttitle\fi\hfill
{\lnum\number\count0}\else\def\\{ and }{\lnum\number\count0}
\hfill\ifx\theshortauthors\relax 
\theauthors\else\theshortauthors\fi\fi\fi}}\def\@evenhead{\@oddhead}
\def\@oddfoot{\small\lfoot\ifnum\count0=\startpage\copyright\ \gtp\hfill\else
\agt, Volume \thevolumenumber\ (\thevolumeyear)\hfill\fi}
\def\@evenfoot{\@oddfoot}
\makeatother
\fi
\let\maketitlepage\makeagttitle
\let\makeshorttitle\maketitlepage
\let\maketitle\maketitlepage

\newwrite\gtoutfile
\long\gdef\makeheadfile{  
{\def\\{, }\def\s{ }
\immediate\openout\gtoutfile head.xxx
\immediate\write\gtoutfile{To: math@arxiv.org}
\immediate\write\gtoutfile{Subject: put OR rep NNNNN:ppppp}
\immediate\write\gtoutfile{--text follows this line--}
\immediate\write\gtoutfile{Proxy-for: \ifx\theasciiauthors\relax
\theauthors\else\theasciiauthors\fi\s<\ifx\theasciiemail\relax\theemail\else\theasciiemail\fi>}
\immediate\write\gtoutfile{\noexpand\\}
\immediate\write\gtoutfile{Authors: \ifx\theasciiauthors\relax
\theauthors\else\theasciiauthors\fi}
{\def\\{ }\immediate\write\gtoutfile{Title: \ifx\theasciititle\relax
\thetitle\else\theasciititle\fi}}
\immediate\write\gtoutfile{Subj-class: GT or SG, GR etc}
\immediate\write\gtoutfile{MSC-class: \theprimaryclass\ifx\thesecondaryclass\relax\else, \thesecondaryclass\fi}
\immediate\write\gtoutfile{Journal-ref: Algebr. Geom. Topol. \thevolumenumber\s
(\thevolumeyear) \startpage-\finishpage}
\immediate\write\gtoutfile{Comments: Published by Algebraic and
Geometric Topology at}
\immediate\write\gtoutfile{\s\s\s  http://www.maths.warwick.ac.uk/agt/AGTVol\thevolumenumber/agt-\thevolumenumber-\thepapernumber.abs.html}
\immediate\write\gtoutfile{\noexpand\\}
\immediate\write\gtoutfile{}
\ifx\theasciiabstract\relax
\immediate\write\gtoutfile{\theabstract}\else
\immediate\write\gtoutfile{\theasciiabstract}\fi
\immediate\write\gtoutfile{}
\immediate\write\gtoutfile{\noexpand\\}
\immediate\write\gtoutfile{}
\immediate\closeout\gtoutfile}}  

\def\maketitlepage{\makeagttitle\makeheadfile}
\let\makeshorttitle\maketitlepage
\let\maketitle\maketitlepage

\def\ifplaintex{\expandafter\ifx\csname documentclass\endcsname\relax}

\def\gtp{{\mathsurround=0pt\it $\cal G\mskip-2mu$eometry \&\ 
$\cal T\!\!$opology $\cal P\!$ublications}}  

\def\recd{{\small Received:\qua\receiveddate\ifx\reviseddate\relax
\else\qquad Revised:\qua\reviseddate\fi\par}} 

\def\lognumber#1{\def\thelognumber{#1}}
\def\volumenumber#1{\def\thevolumenumber{#1}}
\def\volumeyear#1{\def\thevolumeyear{#1}}
\def\papernumber#1{\def\thepapernumber{#1}}
\def\pagenumbers#1#2{\def\startpage{#1}\def\finishpage{#2}}
\def\published#1{\def\publishdate{#1}}

\def\received#1{\def\receiveddate{#1}}

\def\accepted#1{\def\accepteddate{#1}}
\def\asciititle#1{\def\theasciititle{#1}}
\def\covertitle#1{\def\thecovertitle{#1}}

\def\asciiaddress#1{\def\theasciiaddress{#1}}

\long\def\asciiabstract#1{\long\def\theasciiabstract{#1}}
\def\asciikeywords#1{\def\theasciikeywords{#1}}

\let\\\par\let\thelognumber\relax\let\thevolumenumber\relax
\let\thepapernumber\relax\let\thevolumeyear\relax\let\startpage\relax
\let\finishpage\relax\let\publishdate\relax\let\receiveddate\relax
\let\reviseddate\relax\let\accepteddate\relax\let\theasciititle\relax
\let\thecovertitle\relax\let\theasciiauthors\relax\let\theasciiaddress\relax
\let\theasciiabstract\relax\let\theasciikeywords\relax

\let\theasciiemail\relax

\ifplaintex
\font\logobig=cmssbx10 scaled 3836
\font\logomed=cmssbx10 scaled 2557
\else
\font\logobig=cmssbx10 scaled 4200
\font\logomed=cmssbx10 scaled 2800
\fi

\long\def\makeagttitle{   
\count0=\startpage
\agt\hfill      
\hbox to 45truept{\vbox to 0pt{\vglue -13truept{\logomed A\kern -.37em{\logobig 
T}\kern -.38em G}\vss}\hss}
\break
{\small Volume \thevolumenumber\ (\thevolumeyear)
\startpage--\finishpage\nl
Published: \publishdate}

\vglue .25truein

{\parskip=0pt\leftskip 0pt plus
1fil\def\\{\par\smallskip}{\Large\bf\thetitle}\par\medskip} \vglue
0.05truein

{\parskip=0pt\leftskip 0pt plus 1fil\def\\{\par}{\sc\theauthors}
\par\medskip}%
 
\vglue 0.03truein 

{\small\leftskip 25truept\rightskip 25truept{\bf Abstract}\stdspace\theabstract

{\bf AMS Classification}\stdspace\theprimaryclass
\ifx\thesecondaryclass\relax\else; \thesecondaryclass\fi\par
{\bf Keywords}\stdspace \thekeywords\par}\vglue 7truept

}   

\ifplaintex
\font\phead=cmsl9 scaled 950
\font\pnum=cmbx10 scaled 913
\font\pfoot=cmsl9 scaled 950
\headline{\vbox to 0pt{\vskip -4.5mm\line{\small\phead\ifnum
\count0=\startpage ISSN 1472-2739 (on-line) 1472-2747 (printed)
\hfill {\pnum\folio}\else\ifodd\count0\def\\{ }%
\ifx\theshorttitle\relax\thetitle\else\theshorttitle\fi\hfill{\pnum\folio}
\else\def\\{ and }{\pnum\folio}\hfill\ifx\theshortauthors\relax\theauthors
\else\theshortauthors\fi\fi\fi}\vss}}
\footline{\vbox to 0pt{\vglue 0mm\line{\small\pfoot\ifnum\count0=\startpage
\copyright\ \gtp\hfill\else
\agt, Volume \thevolumenumber\ (\thevolumeyear)\hfill\fi}\vss}}
\else
\font=cmsl9 scaled 1050
\font\lnum=cmbx10 
\font=cmsl9 scaled 1050
\makeatletter
\def\@oddhead{{\small\lhead\ifnum\count0=\startpage ISSN 1472-2739 
(on-line) 1472-2747 (printed)\hfill {\lnum\number\count0}\else\ifodd\count0
\def\\{ }\ifx\theshorttitle\relax \thetitle \else\theshorttitle\fi\hfill
{\lnum\number\count0}\else\def\\{ and }{\lnum\number\count0}
\hfill\ifx\theshortauthors\relax 
\theauthors\else\theshortauthors\fi\fi\fi}}\def\@evenhead{\@oddhead}
\def\@oddfoot{\small\lfoot\ifnum\count0=\startpage\copyright\ \gtp\hfill\else
\agt, Volume \thevolumenumber\ (\thevolumeyear)\hfill\fi}
\def\@evenfoot{\@oddfoot}
\makeatother
\fi
\let\maketitlepage\makeagttitle
\let\makeshorttitle\maketitlepage
\let\maketitle\maketitlepage

\newwrite\gtoutfile
\long\gdef\makeheadfile{  
{\def\\{, }\def\s{ }
\immediate\openout\gtoutfile head.xxx
\immediate\write\gtoutfile{To: math@arxiv.org}
\immediate\write\gtoutfile{Subject: put OR rep NNNNN:ppppp}
\immediate\write\gtoutfile{--text follows this line--}
\immediate\write\gtoutfile{Proxy-for: \ifx\theasciiauthors\relax
\theauthors\else\theasciiauthors\fi\s<\ifx\theasciiemail\relax\theemail\else\theasciiemail\fi>}
\immediate\write\gtoutfile{\noexpand\\}
\immediate\write\gtoutfile{Authors: \ifx\theasciiauthors\relax
\theauthors\else\theasciiauthors\fi}
{\def\\{ }\immediate\write\gtoutfile{Title: \ifx\theasciititle\relax
\thetitle\else\theasciititle\fi}}
\immediate\write\gtoutfile{Subj-class: GT or SG, GR etc}
\immediate\write\gtoutfile{MSC-class: \theprimaryclass\ifx\thesecondaryclass\relax\else, \thesecondaryclass\fi}
\immediate\write\gtoutfile{Journal-ref: Algebr. Geom. Topol. \thevolumenumber\s
(\thevolumeyear) \startpage-\finishpage}
\immediate\write\gtoutfile{Comments: Published by Algebraic and
Geometric Topology at}
\immediate\write\gtoutfile{\s\s\s  http://www.maths.warwick.ac.uk/agt/AGTVol\thevolumenumber/agt-\thevolumenumber-\thepapernumber.abs.html}
\immediate\write\gtoutfile{\noexpand\\}
\immediate\write\gtoutfile{}
\ifx\theasciiabstract\relax
\immediate\write\gtoutfile{\theabstract}\else
\immediate\write\gtoutfile{\theasciiabstract}\fi
\immediate\write\gtoutfile{}
\immediate\write\gtoutfile{\noexpand\\}
\immediate\write\gtoutfile{}
\immediate\closeout\gtoutfile}}  

\def\maketitlepage{\makeagttitle\makeheadfile}
\let\makeshorttitle\maketitlepage
\let\maketitle\maketitlepage

\lognumber{5}
\volumenumber{2}
\volumeyear{2002}
\papernumber{5}
\published{5 February 2002}
\pagenumbers{51}{93}
\received{17 October 2001}
\accepted{1 February 2002}

\usepackage{diagrams,amsmath,amssymb}

\newtheorem{theorem}{Th\'eor\`eme}[section]
\newtheorem{proposition}{Proposition}[section]
\newtheorem{coroll}{Corollaire}[section]
\newtheorem{lemma}{Lemme}[section]
\theoremstyle{remark}

\newtheorem{definition}{D\'efinition}[section]

\newcommand{\ope}{\mathcal O} 
\newcommand{\zed}{\mathbb Z}

\newcommand{\Mdg}{\mathbf {R\mathcal M_{dg}}}

\begin{document}

\title {Formes diff\'erentielles g\'en\'eralis\'ees sur une op\'erade\\et
mod\`eles alg\'ebriques des fibrations}
\shorttitle{Formes diff\'erentielles g\'en\'eralis\'ees sur une op\'erade}

\covertitle {Formes diff\noexpand\'erentielles g\noexpand\'en\noexpand\'eralis\noexpand\'ees sur une op\noexpand\'erade\\ 
mod\noexpand\`eles alg\noexpand\'ebriques des fibrations}
\asciititle {Formes differentielles generalisees sur une operade 
et modeles algebriques des fibrations}

\author{David Chataur}
\address{Centre di Recerca Matem\'atica, Institut d'Estudis 
Catalans\\Apartat 50 E - 08193 Bellaterra, Spain}
\asciiaddress{Centre di Recerca Matematica, Institut d'Estudis 
Catalans\\Apartat 50 E - 08193 Bellaterra, Spain}
\email{dchataur@crm.es}

\begin{abstract}
We construct functors of generalized differential forms. In the case
of nilpotent spaces of finite type, they determine the weak homotopy
type of the spaces.  Moreover they are equipped, in an elementary and
natural way, with the action of cup-i products.  Working with
commutative algebras up to homotopy (viewed as algebras over a
cofibrant resolution of the operad of commutative algebras), we show
using these functors that the model of the fiber of a simplicial map
is the cofiber of the algebraic model of this map.

{\bf Resum\'e}\qua On construit des foncteurs de formes
diff\'erentielles g\'en\'eralis\'ees. Ceux-ci, dans le cas d'espaces
nilpotents de type fini, d\'eterminent le type d'homotopie faible des
espaces. Ils sont munis, d'une mani\`ere \'el\'ementaire et naturelle,
de l'action de cup-i produits. Pour les alg\`ebres commutatives \`a
homotopit pr\'es (alg\`ebres sur une r\'esolution cofibrante de
l'op\'erade des alg\`ebres commutatives), on d\'emontre en utilisant
les formes diff\'erentielles g\'en\'eralis\'ees que le mod\`ele de la
fibre d'une application simpliciale est la cofibre du mod\`ele de ce
morphisme.
\end{abstract}

\asciiabstract{On construit des foncteurs de formes differentielles
generalisees. Ceux-ci, dans le cas d'espaces nilpotents de type
fini, determinent le type d'homotopie faible des espaces. Ils sont
munis, d'une maniere elementaire et naturelle, de l'action de
cup-i produits. Pour les algebres commutatives a homotopit pres
(algebres sur une resolution cofibrante de l'operade des
algebres commutatives), on demontre en utilisant les formes
differentielles generalisees que le modele de la fibre d'une
application simpliciale est la cofibre du modele de ce morphisme.

We construct functors of generalized differential forms. In the case
of nilpotent spaces of finite type, they determine the weak homotopy
type of the spaces.  Moreover they are equipped, in an elementary and
natural way, with the action of cup-i products.  Working with
commutative algebras up to homotopy (viewed as algebras over a
cofibrant resolution of the operad of commutative algebras), we show
using these functors that the model of the fiber of a simplicial map
is the cofiber of the algebraic model of this map.}


\asciikeywords{Modeles algebriques, formes differentielles, operades, suites spectrales}

\maketitle


\begin{section}{Introduction}

Il est bien connu que les cocha\^{\i}nes singuli\`eres  d'un espace
topologique forment une alg\`ebre diff\'erentielle gradu\'ee. Celle-ci n'est pas
commutative, mais commutative \`a homotopie pr\`es, ce d\'efaut de commutativit\'e se traduit par l'existence
d'op\'erations sur la cohomologie \`a coefficients dans $\mathbb F_p$: les op\'erations de Steenrod.

En outre, la th\'eorie des op\'erades fournit un cadre agr\'eable pour la construction de mod\`eles alg\'ebriques
pour les espaces topologiques.
Une op\'erade est une structure alg\'ebrique permettant de coder un type
d'alg\`ebres (en particulier des types d'alg\`ebres diff\'erentielles $\mathbb Z$-gradu\'ees). 
La th\'eorie des alg\`ebres sur une op\'erade semble tout adapt\'ee \`a
l'\'etude des structures alg\'ebriques \`a homotopie pr\`es. Pour les alg\`ebres commutatives 
(alg\`ebres sur l'op\'erade
$\mathcal Com$) on a la notion d'alg\`ebres commutatives \`a homotopie pr\`es ou
$E_{\infty}$-alg\`ebres \cite {KM}.

Aussi, V. Hinich et V. Schechtman ont montr\'e que l'ont pouvait munir les 
cocha\^{\i}nes singuli\`eres d'un espace topologique, de mani\`ere naturelle,
d'une structure de $E_{\infty}$-alg\`ebre {\cite {HS}}.

L'homotopie des $E_{\infty}$-alg\`ebres et leur structure de cat\'egorie mod\`ele ferm\'ee
ont \'et\'e \'etudi\'ees r\'ecemment par M. Mandell {\cite {M1}} et V. Hinich {\cite {Hi}}, {\cite {Hi2}}. 
V. Hinich a prouv\'e que la cat\'egorie des op\'erades est une cat\'egorie mod\`ele ferm\'ee, et
qu'une telle structure existe pour les alg\`ebres sur une op\'erade cofibrante. On fixe
un remplacement cofibrant $\mathcal E_{\infty}$ de l'op\'erade des alg\`ebres commutatives; il est alors possible 
d'\'etudier l'homotopie des $\mathcal E_{\infty}$-alg\`ebres, c'est-\`a dire l'homotopie des
alg\`ebres commutatives \`a homotopie pr\`es. 

Pour une op\'erade $E_{\infty}$ particuli\`ere $\mathcal C$ (obtenue \`a partir de l'op\'erade
 des cha\^{\i}nes singuli\`eres de 
l'op\'erade des isom\'etries lin\'eaires), 
M. Mandell a montr\'e 
une \'equivalence de cat\'egories entre une sous-cat\'egorie pleine de  
la cat\'egorie homotopique des $\mathcal C$-alg\`ebres sur $\overline {\mathbb F_p}$ et
la cat\'egorie homotopique des espaces topol\-og\-iques nilpotents $\overline {\mathbb F_p}$-complets. 
V. Hinich a \'etendu ce r\'esultat au cadre des $\mathcal E_{\infty}$-alg\`ebres \cite {Hi2}.
Enfin, M. Mandell a prouv\'e que deux espaces nilpotents de type finis sont faiblement homotopiquement
\'equivalents si et seulement si leurs alg\`ebres de cocha\^{\i}nes singuli\`eres sont 
$\mathcal E_{\infty}$-quasi-isomorphes \cite {M2}.

{\bf R\'esultats}\qua
Le cadre dans lequel on travaille est celui des op\'erades unitales
augment\'ees dans la cat\'egorie des $R$-modules diff\'erentiels $\mathbb Z$-gradu\'es.
On construit des foncteurs $\Omega^{\mathcal O}$ de {\sl
formes diff\'erentielles  g\'en\'eralis\'ees} pour les alg\`ebres sur une
op\'erade $\mathcal O$ (chapitre $3$). Ces foncteurs sont l'analogue  pour les
 $\mathcal O$-alg\`ebres du foncteur $\mathcal A_{Pl}$ de Sullivan pour les
alg\`ebres diff\'erentielles gradu\'ees commutatives. 
On montre que sous des hypoth\`eses raisonnables
(essentiellement que l'op\'erade choisie soit cofibrante) ces foncteurs permettent de
calculer la cohomologie singuli\`ere \`a coefficients (th\'eor\`eme $3.1$):

\begin{theorem}
Pour tout ensemble simplicial $X$ on a $\pi^*\Omega^{\mathcal O}(X)\cong H^*(X,R)$.
\end{theorem} 

La preuve de ce r\'esultat repose essentiellement sur des techniques de mod\`eles 
acycliques adapt\'ees au cadre des $R$-modules diff\'erentiels $\mathbb Z$-gradu\'es.
Toujours en utilisant la th\'eorie des mod\`eles acycliques et la structure de cat\'egorie mod\`ele ferm\'ee
pour les op\'erades unitales augment\'ees, on montre que les foncteurs {\sl
formes diff\'erentielles g\'en\'eralis\'ees}  sont tous munis d'une structure 
de $\mathcal E_{\infty}$-alg\`ebre.
Ce r\'esultat g\'en\'eralise ceux de V. Hinich et V. Schechtman obtenus pour les
cocha\^{\i}nes singuli\`eres  \cite {HS}.
Si on travaille avec des $\mathcal O$-alg\`ebres \`a coefficients dans $\mathbb F$ un corps de
caract\'eristique positive,
cette structure de $\mathcal E_{\infty}$-alg\`ebre implique le r\'esultat suivant (th\'eor\`eme 3.4):

\begin{theorem}
Pour tout ensemble simplicial $X$ l'alg\`ebre $\Omega^{\mathcal O}(X)$ est munie de l'action
d'op\'erations. Cette action induit sur $\pi^*\Omega^{\mathcal O}(X)$ une structure d'alg\`ebre
instable sur l'alg\`ebre de Steenrod telle que  $\pi^*\Omega^{\mathcal O}(X)\cong
H^*(X,\mathbb F)$ soit
un isomorphisme d'alg\`ebres instables.
\end{theorem}

Un des int\'er\^ets des {\sl formes diff\'erentielles g\'en\'eralis\'ees} est qu'elle permettent
de travailler avec n'importe quel remplacement cofibrant de l'op\'erade $\mathcal Com$.
De plus la structure de $\mathcal E_{\infty}$-alg\`ebre n'est pas tr\`es ``lisible" sur les cocha\^{\i}nes
singuli\`eres (les r\'esultats de V. Hinich et V. Schechtman donnent seulement son existence, pour une
description combinatoire de cette structure on pourra se r\'ef\'erer aux travaux de C. Berger et B. 
Fresse \cite {BF}); 
alors que celle-ci est imm\'ediate pour le foncteur des formes diff\'erentielles g\'en\'eralis\'ees
en  $\mathcal E_{\infty}$-alg\`ebres. Pour un choix judicieux d'une op\'erade $E_{\infty}$, on peut
donner des repr\'esentants canoniques pour les cup-i produits.

Ces formes diff\'erentielles permettent de construire, de mani\`ere \'el\'ementaire, une paire de foncteurs
adjoints (en fait une paire de foncteurs adjoints de Quillen) entre la cat\'egorie
des ensembles simpliciaux et de nombreuses cat\'egories d'alg\-\`ebres sur une op\'erade.
Enfin, comme en th\'eorie de Sullivan, le foncteur $\Omega^{\mathcal O}$ permet de d\'efinir un objet
chemin naturel pour les $\mathcal O$-alg\`ebres, et aussi 
des espaces fonctionnels simpliciaux pour ces m\^emes $\mathcal O$-alg\`ebres.

Gr\^ace \`a l'introduction de formes diff\'erentielles g\'en\'eralis\'ees \`a coefficients locaux et
\`a l'extension de ces th\'eories aux ensembles bisimpliciaux, on donne une construction
de la suite spectrale de Leray-Serre (th\'eor\`eme 4.1):

\begin{theorem}
Soit $f:E\longrightarrow B$ une fibration entre ensembles simpliciaux de fibre $F$ (On suppose que la
cohomologie de la base ou celle de la fibre est finie).
Il existe une suite spectrale qui
converge vers $H^*(E;k)$ telle que  $$E^{r,s}_0=\Omega^{\mathcal
O,r}(B;{\mathcal F^s})$$  pour un certain syst\`eme de coefficients locaux
$\mathcal F$. Quand celui-ci est simple (par exemple si $B$ est $1$-connexe)
le terme $E_2$ de cette suite spectrale s'\'ecrit: $$E^{r,s}_2=H^r(B;H^s(F;k)).$$  \end{theorem}

Pour les $\mathcal E_{\infty}$-alg\`ebres cette suite spectrale nous permet de construire un
mod\`ele alg\'ebrique de la fibre d'une application simpliciale cf th\'eor\`eme $4.2$. Ainsi, on obtient 
un r\'esultat qui est classique en homotopie rationnelle: Le mod\`ele alg\'ebrique de la fibre correspond
\`a la cofibre du mod\`ele.
Le th\'eor\`eme $4.2$ s'inscrit dans la lign\'ee des travaux de K. Hess et N. Dupont \cite {DH}.
Il poss\`ede l'avantage de rester valide sur les entiers, mais aussi de donner un mod\`ele de la fibre 
qui d\'etermine son type d'homotopie. On esp\`ere aussi appliquer ce r\'esultat \`a l'\'etude des espaces 
de lacets et des espaces de lacets libres.  
On retrouve aussi dans ce cadre alg\'ebrique les r\'esultats de Kudo sur la transgression dans la suite
spectrale de Leray-Serre.

{\bf Remerciements}\qua Ce travail est en grande partie issu de ma th\`ese ``Formes
diff\'erentielles sur une op\'erade
et mod\`eles alg\'ebriques pour les espaces topol\-og\-iques". J'en profite donc pour remercier
Jean-Louis Loday et Lionel Schwartz dont les conseils de r\'edaction m'ont \'et\'e tr\`es utiles. Je tiens
aussi \`a exprimer toute ma gratitude \`a Beno\^{\i}t Fresse pour ses nombreux conseils et l'attention bienveillante
qu'il apporte \`a mes recherches. Enfin je remercie tr\`es chaleureusement mon directeur de th\`ese Marc Aubry
pour sa disponibilit\'e, son soutien infaillible et ses encouragements constants.
\end{section}

\begin{section}{Homotopie des op\'erades et des alg\`ebres sur une op\'erade}

\begin{subsection}{Op\'erades}

La structure d'{\it op\'erade} a \'et\'e utilis\'ee et introduite au d\'ebut des ann\'ees 1970
par J. M. Boardman, R. M. Vogt \cite {BV} et J. P. May \cite {May1}
dans un contexte
topologique afin d'\'etudier l'homotopie des espaces de lacets it\'er\'es. V. Ginzburg, M. Kapranov \cite {GK} et 
E. Getzler, J.D.S. Jones \cite {Ge} ont repris cette notion dans le cadre alg\'ebrique. On peut \'egalement se r\'ef\'erer
\`a l'expos\'e Bourbaki de J.L. Loday \cite {Lo} ou \`a la premi\`ere partie du livre de I. Kriz et J. P. May
\cite {KM}.

{\bf \noindent La cat\'egorie des $R$-modules diff\'erentiels $\zed$-gradu\'es}\qua 
On fixe un anneau $R$, que l'on suppose commutatif et unitaire.
Les op\'erades avec lesquelles on travaille sont des op\'erades dans la {\it cat\'egorie mono\"{\i}dale 
sym\'etrique des $R$-modules diff\'erentiels $\zed$-gradu\'es} que l'on note $\Mdg$. Par convention,
la diff\'erentielle augmente le degr\'e de 1. 

Soit $M$ un objet de la cat\'egorie $\Mdg$. On note $M^n$ avec $n\in \zed$ le $R$-module des cocha\^{\i}nes de degr\'e 
$n$.

On forme un complexe
$$\ldots\stackrel{d}{\longrightarrow}M^{n-1}\stackrel{d}{\longrightarrow}M^n\stackrel{d}{\longrightarrow}\ldots$$
avec les composantes et la diff\'erentielle de notre $R$-module. Le $R$-module gradu\'e $\pi^*M$ est la cohomologie 
de ce complexe.

Le produit tensoriel de deux $R$-modules diff\'erentiels 
$\zed$-gradu\'es $M$ et $N$ s'\'ecrit:
$$(M\otimes N)^n=\bigoplus_{p+q=n}M^p\otimes N^q.$$
La diff\'erentielle d'un tenseur \'etant donn\'ee par la formule: 
$$d(m\otimes n)=dm\otimes n+(-1)^{|m|}m\otimes dn.$$
\noindent On a un op\'erateur d'\'echange:
$$T:A\otimes B\longrightarrow B \otimes A$$
$$a\otimes b\mapsto (-1)^{|a||b|}b\otimes a.$$
Si $A$ est un $R$-module diff\'erentiel gradu\'e, alors $A[m]$, $m\in \mathbb Z$, d\'esigne la $m$-suspension
de $A$, le $R$-module diff\'erentiel gradu\'e tel que $A[m]^*=A^{m+*}$.

{\bf \noindent $\Sigma_*$-modules}\qua Un {\it $\Sigma_*$-module} ou {\it suite sym\'etrique} est la donn\'ee
d'un ensemble $\{M(n)\}_{n\in \mathbb N}$ d'objets de $\Mdg$. Pour tout entier $n$, le $R$-module $M_n$ est
muni d'une action \`a droite du groupe sym\'etrique $\Sigma_n$. On rearque que les $\Sigma_*$-modules forment
une cat\'egorie monoidale sym\'etrique.

{\bf \noindent Op\'erades}\qua
Une op\'erade $\mathcal O$ 
dans la cat\'egorie $\Mdg$ est la donn\'ee d'un $\Sigma_*$-module $\{\mathcal O(n)\}_{n\in \mathbb N}$,
d'un morphisme $\eta :R\longrightarrow {\mathcal O(1)}$ appel\'e morphisme unit\'e 
et de produits de composition 
$$\gamma:{\mathcal O(k)}\otimes {\mathcal O(j_1)}\otimes \ldots \otimes 
{\mathcal O(j_k)}\longrightarrow {\mathcal O(j)}$$
\noindent d\'efinis pour $1\leq k$ et $\sum j_s=j$.
Ces produits de composition sont, dans un certain sens, associatifs, unitaires et \'equivariants
(cf {\cite {KM}, Part 1, section 1}).

{\bf \noindent Op\'erades unitaires augment\'ees}\qua
On travaille avec des op\'erades unitaires augment\'ees. On note $\bf Op_{ua}$ la cat\'egorie form\'ee
par de telles op\'erades.

Un objet de $\bf Op_{ua}$ est une op\'erade $\mathcal O$ avec un morphisme d'augmentation
$$\epsilon:\mathcal O\longrightarrow \mathcal Com.$$
On rappelle que l'op\'erade $\mathcal Com$ est l'op\'erade des alg\`ebres diff\'erentielles gradu\'ees
commutatives. Elle est donn\'ee par la formule $\mathcal Com(i)=R$.

On demande aussi que $\epsilon$ induise un isomorphisme sur les composantes de degr\'es $0$ et $1$:
$$\mathcal O(0)=\mathcal Com(0)=R$$
$$\mathcal O(1)=\mathcal Com(1)=R.$$
Un morphisme de $\bf Op_{ua}$ est un morphisme d'op\'erades
$f:\mathcal O\longrightarrow \mathcal O'$ tel que $f_0=f_1=Id$.

La cat\'egorie $\bf Op_{ua}$ poss\`ede un objet initial qui est aussi un objet final on note
$\mathbb I$ cette op\'erade. Celle-ci est telle que $\mathbb I(0)=\mathbb I(1)=R$ et
$\mathbb I(i)=0$ pour $i>1$. On v\'erifie aussi que cette cat\'egorie est compl\`ete et cocompl\`ete 
(on forme les limites et les colimites dans la cat\'egorie des op\'erades au-dessus de l'op\'erade $\mathcal Com$).

{\noindent\bf Alg\`ebres sur une op\'erade}\qua
On fixe $\mathcal O$ une op\'erade unitale augment\'ee. Une alg\`ebre sur
$\mathcal O$ (on dit aussi $\mathcal O$-alg\`ebre) est un $R$-module
diff\'erentiel $\zed$-gradu\'e $A$ muni de produits d'\'evaluation: $$\theta_j:{\mathcal
O(j)}\otimes A^{\otimes j}\longrightarrow A$$  associatifs, unitaires et
\'equivariants (pour l'action de $\Sigma_j$).

On note $\bf \mathcal O-\mathcal Alg_{dg}$ la cat\'egorie
des $\mathcal O$-alg\`ebres.

Le foncteur oubli de $\bf \mathcal O-\mathcal Alg_{dg}$ vers $\Mdg$ admet un adjoint \`a gauche,
le foncteur $\mathcal O$-alg\`ebre libre que l'on note aussi $\mathcal O$.
Pour tout $R$-module diff\'erentiel $\zed$-gradu\'e $M$, on d\'efinit $\mathcal O$(M) la $\mathcal O$-alg\`ebre libre
sur M par:
$$\mathcal O(M)=\bigoplus_{p\geq 0}\mathcal O(p)\otimes_{R[\Sigma_p]}M^{\otimes p}=R\oplus M
\bigoplus_{p\geq 2}\mathcal O(p)\otimes_{R[\Sigma_p]}M^{\otimes p}.$$ 
Si $\ope$ est une op\'erade dans la cat\'egorie $\Mdg$ alors
$\pi^*\ope$ est une op\'erade dans la cat\'egorie des R-modules $\zed$-gradu\'es. 

De plus, si $A$ est une $\ope$-alg\`ebre alors $\pi^*A$ est une $\pi^*\ope$-alg\`ebre.

\end{subsection}

\begin{subsection}{La th\'eorie homotopique des op\'erades}

Les op\'erades sont un langage pour aborder l'\'etude d'espaces en topologie alg\'eb\-rique.
Il parait naturel de se placer dans un cadre homotopique. Celui de Quillen (\cite {Ds}, \cite {Hov}, \cite {Qu})
semble convenir tout \`a fait.

\noindent {\bf Homotopie des R-modules diff\'erentiels gradu\'es}\qua
D. Quillen a d\'emontr\'e que l'on pouvait munir la cat\'egorie des $R$-modules
diff\'erentiels $\mathbb N$-gradu\'es d'une structure de cat\'egorie mod\`ele ferm\'ee \cite {Qu}. 
Cette construction s'\'etend au cadre des $R$-modules diff\'erentiels $\mathbb Z$-gradu\'es (\cite {Hov},\cite {St}).
On a alors un structure de cat\'egorie mod\`ele pour laquelle les
\'equivalences faibles sont les quasi-isomorphismes et les fibrations sont les surjections. On 
remarque que tous les objets sont fibrants.

\noindent {\bf Sur les R-modules diff\'erentiels gradu\'es cofibrants}\qua
Pour tout $A$ les conditions suivantes sont \'equivalentes (cf \cite {Sp}):

(a)\qua Le complexe $A$ est cofibrant.

(b)\qua Pour tout complexe acyclique $S$ le complexe $Hom(A,S)$ est aussi acyclique.

En cons\'equence, on en d\'eduit que si $A$ est un objet cofibrant, alors pour tout entier $A^n$ est un $R$-module projectif.
R\'eciproquement
tout complexe de $R$-modules projectifs born\'e sup\'erieurement est cofibrant.

\noindent {\bf Homotopie des $\Sigma_*$-modules}\qua
La cat\'egorie des $\Sigma_*$-modules poss\`ede une structure de cat\'egorie mod\`ele ferm\'ee qui provient de
la structure de cat\'egorie mod\`ele ferm\'ee des $R[\Sigma_n]$-modules diff\'erentiels 
gradu\'es. Explicitement:

i)\qua Un morphisme $f_*:M_*\longrightarrow N_*$ est une \'equivalence faible si pour tout
$n\in {\mathbb N}$ l'application $f_n:M_n\longrightarrow N_n$ est un quasi-isomorphisme de 
$R[\Sigma_n]$-module diff\'erentiel gradu\'e.

ii)\qua  Un morphisme $f_*:M_*\longrightarrow N_*$ est une fibration si pour tout $n\in {\mathbb N}$ le morphisme $f_n$ 
est un \'epimorphisme.

iii)\qua  Un morphisme $f_*:M_*\longrightarrow N_*$ est une cofibration si pour tout $n\in {\mathbb N}$ le morphisme $f_n$ 
est une cofibration.

\noindent {\bf Homotopies des op\'erades}

\noindent{\bf Op\'erades unitaires augment\'ees libres}\qua On note ${\bf \Sigma_*mod}^2$ la cat\'egorie des
$\Sigma_*$-modules 2-r\'eduits. Un objet $\mathcal M$ de ${\bf \Sigma_*mod}^2$ est un $\Sigma_*$-module
tel que $\mathcal M_0=\mathcal M_1=0$.

Soit $\mathcal R$ le  $\Sigma_*$-module 2-r\'eduit tel que $\mathcal R_n=R$
pour tout entier $n \geq 2$. On introduit ${\bf \Sigma_*mod}^2_{a}$ la cat\'egorie des
$\Sigma_*$-modules 2-r\'eduits augment\'es comme \'etant la cat\'egorie des
objets de ${\bf \Sigma_*mod}^2$ au-dessus de $\mathcal R$.

On remarque que la cat\'egorie des $\Sigma_*$-modules 2-r\'eduits augment\'es est munie d'une structure de
cat\'egorie mod\`ele ferm\'ee. En effet, ${\bf \Sigma_*mod}^2$
est une cat\'egorie mod\`ele ferm\'ee. De plus on rappelle que si
$X$ est un objet d'une cat\'egorie mod\`ele ferm\'ee $\mathcal C$, alors la cat\'egorie des objets au-dessus
de $X$ est aussi une cat\'egorie mod\`ele ferm\'ee. Les fibrations sont les \'epimorphismes et les 
\'equivalences faibles sont les quasi-isomorphismes.

\noindent
Soit $U:{\bf Op}_{ua}\longrightarrow{\bf \Sigma_*mod}^2_{a}$ le foncteur oubli de la cat\'egorie des
op\'erades unitaires augment\'ees dans la cat\'egorie des $\Sigma_*$-modules 2-r\'eduits augment\'es.
Le foncteur $U$ admet un foncteur adjoint \`a gauche le foncteur op\'erade unitaire augment\'ee libre not\'e 
$\mathbb T_{ua}$ (Appendice B de \cite{BJT}).

Une op\'erade $\mathcal O$ est dite quasi-libre s'il existe un
$\Sigma_*$-module $M_*$ tel que $\mathbb T_{ua}(M_*)$ soit isomorphe \`a
$\mathcal O$ en tant qu'op\'erade gradu\'ee.

V. Hinich montre que les op\'erades forment une cat\'egorie mod\`ele ferm\'ee (cf \cite {Hi}).
Il est possible d'adapter ce r\'esultat au cadre unital augment\'e.

\begin{theorem}
La cat\'egorie des op\'erades unitaires augment\'ees est munie\break d'une structure de cat\'egorie
mod\`ele ferm\'ee pour laquelle 

{\rm i)}\qua  Un morphisme $f:\mathcal O\longrightarrow \mathcal P$ est une \'equivalence faible si pour tout
$n\in {\mathbb N}$ l'application $f(n):\mathcal O(n)\longrightarrow \mathcal P(n)$ est un quasi-isomorphisme de 
$R[\Sigma_n]$-modules diff\'erentiels gradu\'es. 

{\rm ii)}\qua  Un morphisme $f:\mathcal O\longrightarrow \mathcal P$ est une fibration si pour tout 
$n\in {\mathbb N}$ le morphisme $f(n)$ 
est un \'epimorphisme.
\end{theorem}

\begin{proof}
On reprend les arguments de V. Hinich. On transporte la structure de
cat\'egorie mod\`ele ferm\'ee des  $\Sigma_*$-modules $2$-r\'eduits vers les
op\'erades unitaires augment\'ees via le foncteur $\mathbb T_{ua}$.

Il suffit de v\'erifier que si M est un $\Sigma_*$-module
$2$-r\'eduit cofibrant acyclique ($\pi^*M(n)=0$ pour $n \geq 2$) alors pour
toute op\'erade unitaire augment\'ee $\mathcal O$ le morphisme canonique
$$\mathcal O \longrightarrow \mathcal O\amalg \mathbb T_{ua}(M)$$
est un quasi-isomorphisme.

Ce point se d\'emontre facilement par extension
d'une homotopie contractante de $M$ \`a l'op\'erade $\mathbb T_{ua}(M)$.
Les formules d'extension de V. Hinich restent valables dans notre
cadre.
\end{proof}

\noindent On a pour cette structure la caract\'erisation suivante des objets
cofibrants:

\begin{proposition}
Une op\'erade est cofibrante si et seulement si celle-ci est une r\'etraction
d'une op\'erade quasi-libre.
\end{proposition}

\begin{proof}
La proposition repose sur le fait que la cat\'egorie $\mathcal O_{ua}$ est 
engendr\'ee de mani\`ere cofibrante.
\end{proof}

\noindent {\bf Homotopie des alg\`ebres sur une op\'erade}\qua Si on travaille sur un corps $\mathbb K$ 
de caract\'eristique nulle, nous savons que pour toute op\'erade $\mathcal O$ 
la cat\'egorie des $\mathcal O$-alg\`ebres est munie d'une structure de cat\'egorie mod\`ele
ferm\'ee. Les \'equivalences sont les quasi-isomorphismes, les \'epimorphismes sont les
fibrations, on donnera plus loin une description des cofibrations. On appelle cette structure 
la {\sl structure mod\`ele adjointe}. Cette terminologie provient de l'adjonction de Quillen
entre la cat\'egorie des $\mathcal O$-alg\`ebres et la cat\'egorie des $\mathbb K$-espaces vectoriels
 diff\'erentiels gradu\'es.

Si on travaille sur un anneau $R$ quelconque ce r\'esultat n'est pas toujours
vrai. Par exemple, sur $\mathbb F_p$, la cat\'egorie 
des $\mathbb Z$-alg\`ebres diff\'erentielles gradu\'ees commutatives ne peut \^etre
munie d'une structure mod\`ele adjointe (il existe une autre structure de
cat\'egorie mod\`ele ferm\'ee pour ces alg\`ebres \cite {St}). Par contre,
si $\mathcal O$ est une op\'erade; la cat\'egorie des $\mathcal O$-alg\`ebres est
munie d'une structure mod\`ele adjointe si et seulement si pour toute $\mathcal O$-alg\`ebre $A$ 
et tout $n\in\mathbb Z$ le morphisme
canonique 
$$A\longrightarrow A\amalg \mathcal O(x_n,dx_n)$$
est un quasi-isomorphisme.

On retrouve l'analogue du r\'esultat de V. Hinich \cite {Hi} pour les
op\'erades unitaires augment\'ees cofibrantes. En effet, si $\mathcal O$ est une op\'erade cofibrante,
alors la cat\'egorie des $\mathcal O$-alg\`ebres admet une structure mod\`ele adjointe.
Supposons que $\mathcal O$ et $\mathcal O'$ soient deux mod\`eles cofibrants de la m\^eme op\'erade
$\mathbb O$; alors la cat\'egorie des $\mathcal O$-alg\`ebres et la cat\'egorie des 
$\mathcal O'$-alg\`ebres sont \'equivalentes au sens de Quillen.
La preuve repose sur une filtration du coproduit pour une op\'erade
quasi-libre (cf \cite {Hi2}).

{\noindent \bf Extensions libres}\qua
Dans cette  section on caract\'erise les cofibrations pour les $\mathcal O$-alg\`ebres.

Notons $A\amalg B$ le coproduit de deux $\mathcal O$-alg\`ebres (pour une r\'ealisation de
ce coproduit consulter \cite {Fr}). On peut g\'en\'eraliser la
notion d'extension libre (on parle aussi de morphisme quasi-libre)
introduite dans le cadre des alg\`ebres diff\'erentielles gradu\'ees \cite {FHT}
aux  $\mathcal O$-alg\`ebres.

\begin{definition}
Une extension libre est un morphisme de $\mathcal O$-alg\`ebres
$$A\stackrel{i}{\longrightarrow} A\amalg_{\tau}\mathcal O(M)$$
tel que:

a)\qua $A\amalg_{\tau}\mathcal O(M)=A\amalg\mathcal O(M)$ en tant que module gradu\'e.

b)\qua Le morphisme $i$ est l'application canonique.

c)\qua Le $R$-module $M$ s'\'ecrit sous la forme $M=\bigcup^{\infty}_{i=0}M(i)$ avec

$M(i)\subset M(i+1)$.
Les $R$-modules gradu\'es $M(0)$ et $M(i+1)/M(i)$ sont libres.

d)\qua La diff\'erentielle $d$ est telle que $d:M(0)\longrightarrow A$ et $d:M(i+1)\longrightarrow A\amalg
\mathcal O(M(i))$.
\end{definition}

Toute application $f:A\longrightarrow B$ admet une factorisation:
\begin{diagram}
A       &  &\rTo^f &       & B      \\
  &\rdTo^i &              &\ruTo^p &\\
  &      &A \amalg_{\tau} \mathcal O(M) & &       \\
\end{diagram}
avec $i$ une extension libre et $p$ une fibration acyclique.

De plus, un morphisme est une cofibration si et seulement si c'est une r\'etraction d'une extension libre. Ceci repose
sur le fait que la cat\'egorie des $\mathcal O$-alg\`ebres est engendr\'ee de mani\`ere cofibrante.
On peut montrer que la cat\'egorie des $\mathcal O$-alg\`ebres est une cat\'egorie
mod\`ele cellulaire \cite {Hir}, \cite {Hov} (ce qui est aussi le cas de la cat\'egorie des op\'erades unitaires
augment\'ees et
de la cat\'egorie des $\Sigma_*$-modules).

\end{subsection}

\begin{subsection}{Sur les $\mathcal O$-alg\`ebres \`a homotopie pr\`es}

Soit $\mathbb O$ une op\'erade unitaire et augment\'ee, on fixe un mod\`ele
cofibrant  $\mathcal O$ de cette op\'erade. 

\begin{definition}
On appelle $\mathbb O$-alg\`ebre \`a homotopie pr\`es un objet de la
cat\'eg\-orie des $\mathcal O$-alg\`ebres. \end{definition}

\noindent Nous avons vu pr\'ec\'edemment que du point vue homotopique le choix
du mod\`ele cofibrant importait peu (pour deux mod\`eles cofibrants de la
m\^eme op\'erade les cat\'egories homotopiques sont \'equivalentes).

Une des difficult\'es consiste \`a construire un mod\`ele cofibrant explicite pour une op\'erade 
$\mathbb O$ donn\'ee.

Dans le cadre rationnel et pour les op\'erades de Koszul (\cite {GK}, \cite {Lo}), 
il existe un proc\'ed\'e utilisant la construction cobar $B^*$ pour obtenir un tel mod\`ele. 
Si $\mathbb O$ est une op\'erade de Koszul, si $\mathbb O^{!}$ est son dual
de Koszul, alors $B^*(\mathbb O^{!}[-1])^*$ est une r\'esolution quasi-libre
de $\mathbb O$. Toujours dans le cadre rationnel une th\'eorie du mod\`ele
minimal a \'et\'e d\'evelopp\'e par M. Markl \cite {Mar}.

Toute $\mathbb O$-alg\`ebre est une $\mathcal O$-alg\`ebre via le morphisme
d'op\'erades $\mathcal O\longrightarrow \mathbb O$.

On s'int\'eresse plus sp\'ecifiquement au cas des alg\`ebres commutatives \`a homotopie pr\`es.
On choisit un mod\`ele cofibrant de l'op\'erade $\mathcal Com$ que l'on note $\mathcal E_{\infty}$.
On dit que $\mathcal E_{\infty}$ est une $E_{\infty}$-op\'erade cofibrante.
On peut m\^eme supposer que l'op\'erade $\mathcal E_{\infty}$ est quasi-libre.
Dans ce cas, les modules diff\'erentiels gradu\'es $\mathcal E_{\infty}(n)$ sont tous $R[\Sigma_n]$-libres et
acycliques. On peut choisir ces modules diff\'erentiels $\mathbb Z$-gradu\'es tels que 
$\mathcal E_{\infty}(n)^p=0$ pour $p>0$.

Dans le langage de P. May \cite {KM} une $E_{\infty}$-op\'erade $\ope$ est une op\'erade telle que $\ope(l)$ soit
une r\'esolution $R[\Sigma_l]$-projective de $R$.
On remarque qu'une $E_{\infty}$-op\'erade au sens de P. May n'est pas n\'ecessairement
cofibrante en tant qu'op\'erade.

Si on travaille avec une $E_{\infty}$-op\'erade $\ope$ 
et avec $R$ un corps de caract\'eristique positive, alors l'homotopie des $\ope$-alg\`ebres est munie d'op\'erations:

Soient $\ope$ une $E_{\infty}$-op\'erade et $A$ une $\ope$-alg\`ebre, il existe pour tout $s\geq 0$ et 
$R$ de caract\'eristique $2$ des
op\'erations:
$$\mathcal P^s:\pi^qA\longrightarrow \pi^{q+s}A$$
et pour $R$ de caract\'eristique $p>2$:
$$\mathcal P^s:\pi^qA\longrightarrow \pi^{q+2s(p-1)}A.$$
Ces op\'erations v\'erifient les propri\'et\'es suivantes:

i)\qua $\mathcal P^s(x)=0$ si $p=2$ et $s<|x|$ ou si $p>2$ et $2s<|x|$.

ii)\qua $\mathcal P^s(x)=x^p$ si $p=2$ et $s=|x|$ ou si $p>2$ et $2s=|x|$.

iii)\qua $\mathcal P^s(xy)=\sum \mathcal P^t(x)\mathcal P^{s-t}(y)$ (formule de Cartan).

iv) (formule d'Adem) Si $p\geq 2$ et $t>ps$:
$$\mathcal P^t\mathcal P^s=\sum_{i}(-1)^{t+i}(pi-t,t-(p-1)s-i)\mathcal P^{s+t-i-1}\mathcal P^i$$
si $p>2$, $t>ps$, et si par $\beta$ on note le mod-$p$ Bockstein, alors:
$$\mathcal P^t\beta\mathcal P^s=\sum_{i}(-1)^{t+i}(pi-t,t-(p-1)s-i)\beta\mathcal P^{s+t-i-1}\mathcal P^i$$
$$-\sum_i(-1)^{t+i}(pi-t-1,t-(p-1)s-i)\mathcal P^{s+t-i-1}\beta\mathcal P^i$$
$(i,j)=\frac{(i+j)!}{i!j!}$ si $i\geq 0$ et $j\geq 0$ et $(i,j)=0$ si $i$ ou $j$ sont n\'egatifs.

{\noindent \bf Exemple de $E_{\infty}$-op\'erade: la r\'esolution bar des groupes sym\'etriques}\qua Soit $\mathcal R_B$ le $\Sigma_*$-module tel que
$\mathcal R_B(n)$ est le complexe normalis\'e de la r\'esolution bar du groupe sym\'etrique $\Sigma_n$.
Ce $\Sigma_*$-module est une op\'erade. L'op\'erade $\mathcal R_B$ est \'evidemment une
$E_{\infty}$-op\'erade au sens de P. May, mais elle n'est pas cofibrante
(l'op\'erade $\mathcal Com$ est r\'etracte de cette op\'erade). Les
cog\`ebres sur cette  op\'erade ont \'et\'e \'etudi\'ees par J. Smith \cite
{Sm}. V.A. Smirnov a lui ausi \'etudi\'e les  cog\`ebres sur une
$E_{\infty}$-op\'erade et leurs liens avec l'homotopie des espaces
topologiques \cite {Sm85}, \cite {Sm98}. 

\end{subsection}

\end{section}

\begin{section}{Formes diff\'erentielles g\'en\'eralis\'ees}

Dans ce chapitre on construit des foncteurs de formes diff\'erentielles g\'en\'eralis\'ees. Ce sont
des foncteurs de la cat\'egorie des ensembles simpliciaux vers une cat\'eg\-orie de $\mathcal O$-alg\`ebres.

Gr\^ace \`a la th\'eorie de mod\`eles acycliques que nous d\'eveloppons dans le premier
paragraphe on \'etablit une \'equivalence d'homotopie entre ces foncteurs et le foncteur 
des cocha\^{\i}nes singuli\`eres normalis\'ees (c'est le r\'esultat principal du second paragraphe).
 
Dans le troisi\`eme paragraphe, on \'etudie la structure multiplicative de ces foncteurs et on montre
qu'ils sont tous \`a valeurs dans les $\mathcal E_{\infty}$-alg\`ebres.

Enfin dans le dernier paragraphe, on construit une paire de foncteurs adjoints de Quillen entre les ensembles 
simpliciaux et les $\mathcal O$-alg\`ebres, via un foncteur de  formes diff\'erentielles g\'en\'eralis\'ees
et un foncteur de r\'ealisation simpliciale. On donne aussi quelques applications \`a l'homotopie des
$\mathcal O$-alg\`ebres (espace de chemins et homotopie simpliciale). 
 
\begin{subsection}{La th\'eorie des mod\`eles acycliques}
On \'etend la th\'eorie
des mod\`eles acycliques \cite {EM2}, \cite{Ma}, \cite{S} au cadre $\mathbb Z$-gradu\'e.

Soit $F:\mathcal S^{op}\longrightarrow \mathbf C$ un foncteur contravariant des ensembles simpliciaux \`a
valeurs dans une cat\'egorie $\mathbf C$. 
On associe \`a $F$ le foncteur contravariant $F':\mathcal S^{op}\longrightarrow \mathbf C$ tel que:
$$F'(X)=\prod_{x\in X_n}F(\Delta[n])$$
\noindent o\`u $\Delta[n]$ est le simplexe standard de dimension $n$.
Le produit est pris sur tous les $n\geq 0$ et les $x\in X_n$. On notera $\{m_x,x\}\in F'(X)$
l'\'el\'ement dont la composante index\'ee par $x\in X_n$ est l'\'el\'ement $m_x\in F'(\Delta[n])$.
Soit $f:X\longrightarrow Y$
un morphisme entre ensembles simpliciaux.
Le morphisme associ\'e $F'(f):F'(Y)\longrightarrow F'(X)$ est donn\'e par
$$F'(f)\{m_y,y\}=\{m_{f(x),x}\}_{y=f(x)}.$$
\noindent Une transformation naturelle $T:F\longrightarrow G$ induit $T':F'\longrightarrow G'$ donn\'ee par
la formule: 
$$T'(X)\{m_x,x\}=\{T(\Delta[n])m_x,x)\}.$$
\noindent On d\'efinit aussi une transformation naturelle $\Phi:F\longrightarrow F'$ en posant 
$$\Phi(X)u=\{F({\bf x})u,x\}$$ 
pour $u \in F(X)$.
Dans cette formule, on utilise le fait que la donn\'ee de $x\in X_n$ est \'equivalente \`a un
morphisme $\bf x:\Delta[n]\longrightarrow X$. On 
v\'erifie facilement la formule $T'\Phi=\Phi T$.

\begin{definition}
a)\qua On dit que $F:\mathcal S^{op}\longrightarrow \Mdg$ 
est corepr\'esentable s'il existe une transformation naturelle 
$\Psi:F'\longrightarrow F$ telle que $\Psi$ est inverse \`a gauche de $\Phi$.

b)\qua Un foncteur $F:\mathcal S^{op}\longrightarrow \Mdg$ est augment\'e s'il 
existe une transformation naturelle $\epsilon:F\longrightarrow R$. On suppose
que pour tout $n$ le morphisme $\epsilon:F(\Delta[n])\longrightarrow R$ est une fibration.

c)\qua On dit que $F$ est acyclique sur les mod\`eles, si pour tout $n\in \mathbb N$, 
 l'augment\-ation $\epsilon:F(\Delta[n])\longrightarrow R$ est une \'equivalence faible.

d)\qua Le foncteur $F$ est cofibrant, si pour tout $n\in \mathbb N$, 
 $F(\Delta[n])$ est cofibrant.
\end{definition}

\begin{proposition} Supposons que $F:\mathcal S^{op}\longrightarrow \mathbf C$ est un foncteur corepr\'esent\-able,
 augment\'e, acyclique sur les mod\`eles
et que le foncteur $G:\mathcal S^{op}\longrightarrow \mathbf C$ est augment\'e, cofibrant.

Alors, il existe une transformation naturelle $f:G\longrightarrow F$ telle que $\epsilon f=\epsilon$.

De plus, deux transformations naturelles $f,g:G\longrightarrow F$ telles que $\epsilon f=\epsilon =\epsilon g$, sont homotopes.
Plus pr\'ecis\'ement, il existe une homotopie \`a gauche (naturelle) $h:G\longrightarrow F$ entre f et g.
\end{proposition}

\begin{proof}[Preuve]
i)\qua Consid\'erons le diagramme suivant:
\begin{diagram}
 &             & F(\Delta[n])\\
 &\ruDotsto^{f'_n}    & \dTo^{\epsilon} \\
G(\Delta[n])         & \rTo^{\epsilon}  & R\\ 
\end{diagram}
l'application verticale est une fibration acyclique (par hypoth\`ese), et $G(\Delta[n])$ est
cofibrant. Sous ces hypoth\`eses, il existe un rel\`evement $f'_n$.
On obtient ainsi une transformation naturelle
 $f':F'\longrightarrow G'$. On pose $f=\Psi f' \Phi$, c'est la transformation naturelle demand\'ee.

ii)\qua Comme le rel\`evement est unique \`a une homotopie \`a gauche pr\`es, on en d\'eduit une famille 
d'homotopie \`a gauche $H'_n$
de $f'_n$ \`a $g'_n$. 

Comme il existe dans $\Mdg$ un objet cylindre naturel qui commute avec les produits, ceci permet de
construire une  homotopie $H':IF'\longrightarrow G'$ de $f$ \`a $g$. Gr\^ace \`a la naturalit\'e 
de ce m\^eme objet chemin on a une homotopie \`a gauche $H:IF\longrightarrow G$.
\end{proof}

\begin{coroll}
Si $F$ et $G$ sont tous deux corepr\'esentables, augment\'es, cofibrants et
acycliques sur les mod\`eles, alors F et G sont naturellement homotopiquement
\'equivalents.

Toute transformation naturelle entre deux foncteurs corepr\'esentables, augmen\-t\'es, cofibrants et acycliques sur 
les mod\`eles qui commute aux augmentations induit une \'equivalence d'homotopie.
\end{coroll}

On note $\mathcal C:\mathcal S^{op}\longrightarrow \Mdg$ le foncteur des cocha\^{\i}nes singuli\`eres 
normalis\'ees. Le foncteur $\mathcal C$ est corepr\'esentable, augment\'e, cofibrant et acyclique sur les mod\`eles.
On en d\'eduit le r\'esultat suivant:  

\begin{proposition} Soit $F:\mathcal S^{op}\longrightarrow \Mdg$ un foncteur contravariant qui est 
corepr\'esentable, augment\'e, cofibrant et acyclique sur les  mod\`eles, et 
$\mathcal C:\mathcal S^{op}\longrightarrow \Mdg$ 
le foncteur des cocha\^{\i}nes singuli\`eres normalis\'ees.
 
Alors, les foncteurs $F$ et $\mathcal C$ sont naturellement \'equivalents.
Donc, pour tout ensemble simplicial X on a un isomorphisme naturel:
$$\pi^*F(X)\cong \pi^*C(X)\cong H^*(X;R).$$ 
\end{proposition}

\end{subsection}

\begin{subsection}{Le foncteur des formes diff\'erentielles g\'en\'eralis\'ees pour les alg\`ebres sur une op\'erade}

Nous construisons dans ce qui suit un foncteur de formes diff\'erentielles g\'en\'eral\-is\'ees pour les alg\`ebres 
sur une op\'erade $\ope$. On prouve que ce foncteur qui est not\'e $\Omega^{\mathcal O}$ v\'erifie les
hypoth\`eses de la proposition ci-dessus. Pour \^etre plus pr\'ecis, on fixe $\mathcal O$ une op\'erade unitaire
augment\'ee avec les propri\'et\'es suivantes:

a)\qua Le morphisme d'augmentation $\epsilon:\ope\longrightarrow \mathcal Com$ est
une fibration. 

b)\qua Le $R$-module diff\'erentiel $\ope(n)$ est concentr\'e en
degr\'es n\'egatifs. 

c)\qua L'op\'erade $\ope$ est cofibrante.

L'hypoth\`ese a) permet de d\'efinir une suite d'applications $s_n:R\longrightarrow \ope(n)$ telles que
$\epsilon s_n=Id_R$. Les $\{s_n\}_{n\in \mathbb N}$ ne donnent en aucun cas un morphisme d'op\'erades, sinon
l'op\'erade $\mathcal Com$ serait une r\'etraction de $\ope$ et donc cofibrante.

Posons $1_n=s_n(1)$. Comme $\ope(0)=\ope(1)=R$ (car $\mathcal O$ est unitaire augment\'ee), on a
$1_0=1_1=1$.

\begin{proposition}
Soit $A$ une $\ope$-alg\`ebre.
Le morphisme structural $$\theta_2:\ope(2)\otimes A^{\otimes 2}\longrightarrow A$$
 v\'erifie l'\'equation:
$$\theta_2(1_2\otimes\lambda\otimes a)=\theta_2(1_2\otimes a\otimes\lambda)=\lambda .a$$
\noindent pour tout $\lambda \in R$ et tout $a\in A$.
\end{proposition}

\begin{proof}
On consid\`ere le produit $\gamma_2:\ope(2)\otimes \ope(0)\otimes\ope(1)\longrightarrow \ope(1)$.
Le diagramme commutatif:
\begin{diagram}
\ope(2)\otimes \ope(0)\otimes\ope(1) & \rTo^{\gamma_2} & \ope(1)\\
\dTo^{\epsilon\otimes\epsilon\otimes\epsilon}& &\dTo^{\epsilon}\\
\mathcal Com(2)\otimes \mathcal Com(0)\otimes\mathcal Com(1) &\rTo^{\gamma_2} &\mathcal Com(1)\\
\end{diagram}
montre que $\gamma_2(1_2\otimes 1_0\otimes 1_1)=1_1=1$. Les relations de la proposition sont une
cons\'equence de ces identit\'es.
\end{proof}

\noindent Cette proposition montre que pour toute $\ope$-alg\`ebre $A$
le produit:
$$\mu:A\otimes A\longrightarrow A$$
$$\mu(a\otimes b)=\theta_2(1_2\otimes a \otimes b)$$
poss\`ede une unit\'e.
Cette propri\'et\'e d'unitalit\'e n'est pas v\'erifi\'ee pour une
$E_{\infty}$-op\'erade quelconque.

\begin{definition}
On a une $\ope$-alg\`ebre simpliciale ${\Omega^{\ope}}_*$.
dont la composante de dimension simpliciale $n$ est la $\ope$-alg\`ebre
$${\Omega^{\ope}}_n=\frac {\ope(x_0,\ldots,x_n;dx_0,\ldots,dx_n)}{I_n}$$
engendr\'ee par les \'el\'ements $x_0,\ldots,x_n$ de degr\'e $0$ et
$dx_0,\ldots,dx_n$ de degr\'e $1$ et quotient\'ee par l'id\'eal $I_n$ engendr\'e par les
relations $\sum x_i=1_0$, $\sum dx_i=0$.
Les op\'erateurs de face sont donn\'es par les formules:

$\delta_j(x_j)=0$ et $\delta_j(dx_j)=0$;

$\delta_j(x_i)=x_i$ et $\delta_j(dx_i)=dx_i$ si $i<j$;

$\delta_j(x_i)=x_{i-1}$ et $\delta_j(dx_i)=dx_{i-1}$ si $i>j$;

les op\'erateurs de d\'eg\'en\'erescence par les formules:

$\sigma_j(x_j)=x_j+x_{j+1}$ et $\sigma_j(dx_j)=dx_j+dx_{j+1}$;

$\sigma_j(x_i)=x_i$ et $\sigma_j(dx_i)=dx_i$ si $i<j$;

$\sigma_j(x_i)=x_{i+1}$ et $\sigma_j(dx_i)=dx_{i+1}$ si $i>j$.

On d\'efinit le foncteur des formes diff\'erentielles g\'en\'eralis\'ees:
$${\Omega^{\ope}}:\mathcal S^{op} \longrightarrow \mathcal O-Alg_{dg}$$
par la formule:
$${\Omega^{\ope}}(X)=Hom_{\mathcal S}(X,{\Omega^{\ope}}_*)$$
\end{definition}

On remarque que ${\Omega^{\ope}}(\Delta[n])=Hom_{\mathcal S}(\Delta[n],{\Omega^{\ope}}_*)\cong{\Omega^{\ope}}_n$. 
On va montrer que le foncteur des formes 
diff\'erentielles g\'en\'eralis\'ees est corepr\'esentable, cofibrant et acyclique sur les mod\`eles.

Si on applique la th\'eorie des mod\`eles acycliques, on d\'eduit que ce foncteur
est homotopiquement \'equivalent au foncteur des cocha\^{\i}nes singuli\`eres.

\begin{lemma}
Si l'op\'erade $\ope$ est cofibrante et si chaque R-module diff\'erentiel $\zed$-gradu\'e  $\ope (l)$ est 
concentr\'e en degr\'es n\'egatifs, alors ${\Omega^{\ope}}_n$ est cofibrant et acyclique relativement \`a 
l'augmentation.
\end{lemma}

\begin{proof}[Preuve]
Montrons d'abord que ${\Omega^{\ope}}_n$ est acyclique; \`a cette fin on remarque que l'on a 
l'isomorphisme de $\ope$-alg\`ebres:
$$\frac {\ope(x_0,\ldots,x_n;dx_0,\ldots,dx_n)}{I_n}\cong \ope(x_1,\ldots,x_n;dx_1,\ldots,dx_n)$$
Or $\ope(x_1,\ldots,x_n;dx_1,\ldots,dx_n)$ est acyclique. En effet, pour une
op\'erade cofibrante une alg\`ebre libre sur un R-module diff\'erentiel
gradu\'e acyclique est acyclique.

Il reste \`a  montrer que ${\Omega^{\ope}}_n$ est cofibrant
en tant que $R$-module diff\'erentiel $\zed$-gradu\'e.

On rappelle qu'un R-module diff\'erentiel $\zed$-gradu\'e projectif n'est pas
n\'ecessaire\-ment cofibrant \cite {Sp}. Mais c'est le cas s'il est born\'e sup\'erieurement.

En utilisant l'isomorphisme pr\'ec\'edent ${\Omega^{\ope}}_n$, on a:
$${\Omega^{\ope}}_n=\bigoplus_{l\geq0}\ope(l)\otimes_{R[\Sigma_l]}M^{\otimes l}$$
On pose ${\Omega^{\ope}}_n(l)=\ope(l)\otimes_{R[\Sigma_l]}M^{\otimes l}$. Montrons que ce R-module  
diff\'erentiel $\zed$-gradu\'e est cofibrant, on en d\'eduit alors que ${\Omega^{\ope}}_n$ est aussi cofibrant.

Comme $\ope$ est cofibrante il existe une op\'erade quasi-libre $\ope^"$ telle que $\ope$ est une r\'etraction de
$\ope"$. Donc pour tout $l$, $\ope(l)$ est un r\'etract de $\ope^"(l)$. Et $\ope^"(l)$ est un 
$R[\Sigma_l]$-module libre, car $\ope^"$ est une op\'erade libre et cofibrante. 
On peut supposer que cette op\'erade est aussi born\'ee sup\'erieurement. 
Alors ${\Omega^{\ope}}_n(l)$ est un r\'etract de $\ope^"(l)\otimes_{R[\Sigma_l]}M^{\otimes l}$. L'objet $M$
est un R-module diff\'erentiel $\zed$-gradu\'e libre concentr\'e en degr\'es 0 et 1.
On en d\'eduit facilement que
$\ope^"(l)\otimes_{R[\Sigma_l]}M^{\otimes l}$ est cofibrant (car celui-ci est $R$-libre 
en tout degr\'e et born\'e sup\'erieurement). 

Le r\'esultat est une cons\'equence imm\'ediate du fait que les objets cofibrants
sont stables par r\'etraction.
\end{proof}

Afin de montrer que le foncteur des formes diff\'erentielles g\'en\'eralis\'ees est corep\-r\'es\-entable,
il suffit de prouver que l'alg\`ebre simpliciale ${\Omega^{\ope}}_*$ est contractile en tant qu'ensemble
simplicial. En effet d'apr\`es M. Majewski \cite {Ma}, si $\mathbb A:\mathcal S\longrightarrow \Mdg$ est un foncteur 
de la cat\'egorie
des ensembles simpliciaux \`a valeurs dans $\Mdg$ tel que $\mathbb A(X)=Hom_{\mathcal S}(X,M^*)$, alors $\mathbb A$
est corepr\'esentable si $M^*$ est contractile.

\begin{lemma} Pour tout $s\geq 1$, les groupes d'homotopie $\pi_{s-1}({\Omega^{\ope}}_*)$ sont triviaux.
\end{lemma}

\begin{proof}[Preuve]On \'etend \`a l'alg\`ebre simpliciale ${\Omega^{\ope}}_*$
la preuve donn\'ee par M. Karoubi dans le cadre des 
formes diff\'erentielles non commutatives (\cite{Ka1},\cite{Ka2}). 
On commence par montrer que  $\pi_0({\Omega^{\ope}}_*)$ est trivial. Soit $\omega \in {\Omega^{\ope}}_0$. On
identifie ${\Omega^{\ope}}_0$ avec $R$, on suppose que $\omega$ est un scalaire. On identifie aussi
${\Omega^{\ope}}_1$ avec $\mathcal O(x,dx)$.
L'application $\delta_0$ correspond \`a l'\'evaluation en $x=0$, $dx=0$ et $\delta_1$
\`a l'\'evaluation $x=1$ et $dx=0$.
L'\'el\'ement $\theta=\omega (1-x)$ v\'erifie $\delta_0\theta=\omega$, $\delta_1\theta=\omega$. D'o\`u le
r\'esultat.
On suppose $s\geq 2$. Montrer que $\pi_{s-1}({\Omega^{\ope}}_*)$ est trivial est
\'equivalent \`a prouver que, pour toute forme $\omega\in {\Omega^{\ope}}_{s-1}$ satisfaisant $\delta_i\omega=0$ 
pour tout $i$, il
existe $\theta\in {\Omega^{\ope}}_s$ tel que $\delta_0\theta=\omega$ et $\delta_i\theta=0$ si $i>0$.
Comme on a $1-t_1-\ldots -t_s=0$ et $dt_1+\ldots +dt_s=0$ dans ${\Omega^{\ope}}_{s-1}$ on remplace $t_1$
par $1-t_2-\ldots -t_s$ et $dt_1$ par $-dt_2-\ldots -dt_s$. Puis on pose
$$\omega =\mu_1 (t_2,\ldots ,t_s;dt_2,\ldots ,dt_s)$$ 
avec $\mu_1 \in \ope (t_2,\ldots ,t_s;dt_2,\ldots ,dt_s)$.
Enfin on d\'efinit $\sigma\in {\Omega^{\ope}}_{s}$:
$$\sigma_1 (t_0,t_1,\ldots ,t_s;dt_0,\ldots ,dt_s)=
\gamma_2(1_2\otimes t_1\otimes\mu_1 (t_2,\ldots ,t_s;dt_2,\ldots ,dt_s))$$ 
On v\'erifie que les restrictions aux faces sont nulles pour $i>0$ et que la restriction \`a la $0$-face est de 
la forme $\theta_2(1_2\otimes t_1\otimes\omega (t_1,\ldots,t_s;dt_1,\ldots,dt_s))$.
De mani\`ere analogue on construit pour chaque $i$ des formes $\sigma_i$ telles que:

$\delta_j(\sigma_i)=0$ pour $j\not =0$,

$\delta_0(\sigma_j)=\theta_2(1_2\otimes t_j\otimes\omega (t_1,\ldots ,t_s;dt_1,\ldots ,dt_s))$.
 
Et la forme $\theta$ est donn\'ee par la 
formule suivante:
$$\theta =\sum_{1\leq i \leq s}\sigma_i(t_0,t_1,\ldots ,t_s;dt_0,\ldots ,dt_s).$$ 
Cette forme v\'erifie:

$\delta_j(\theta)=0$ pour $j\not =0$,
$$\delta_0(\theta)=\sum_{1\leq i \leq s}\gamma_2(1_2\otimes t_j\otimes\omega (t_1,\ldots ,t_s;dt_1,\ldots ,dt_s))$$
$$\delta_0(\theta)=\theta_2(1_2\otimes\sum_{1\leq i \leq s} t_j\otimes\omega (t_1,\ldots ,t_s;dt_1,\ldots ,dt_s))$$
$$\delta_0=\theta_2(1_2\otimes 1_0\otimes\omega (t_1,\ldots ,t_s;dt_1,\ldots ,dt_s)).$$
Comme $\theta_2(1_2\otimes 1_0\otimes w)=w$ (d'apr\`es la proposition $3.3$) on a:

$\delta_0(\theta)=\omega (t_1,\ldots ,t_s;dt_1,\ldots ,dt_s).$ 
\end{proof}

\begin{theorem}
Soit $\ope$ une op\'erade cofibrante, unitaire, augment\'ee (l'aug\-men\-tation est surjective),
telle que chaque $R$-module diff\'erentiel gradu\'e $\ope (l)$ est born\'e sup\'erieurement.

Alors l'alg\'ebre ${\Omega^{\ope}}(X)$ est homotopiquement \'equivalente \`a $C^*X$ (comme 
R-module diff\'erentiel  $\zed$-gradu\'e). Cette \'equivalence d'homotopie est naturelle en $X$.
\end{theorem}

\begin{proof}[Preuve]
C'est une cons\'equence imm\'ediate des lemmes pr\'ec\'edents et de la th\'eorie des mod\`eles acycliques.
\end{proof}

\noindent{\bf Remarques}\qua On peut en fait montrer que si $\ope$ est une op\'erade unitaire, augment\'ee 
(l'augmentation est toujours suppos\'ee surjective), telle que chaque $R$-module 
diff\'erentiel gradu\'e $\ope (l)$ 
est born\'e sup\'erieurement et $R[\Sigma_l]$-projectif, alors le th\'eor\`eme pr\'ec\'edent est encore valable.
En particulier cette construction s'applique \`a l'op\'erade $\mathcal R_B$.
Cette construction s'applique aussi \`a $\mathcal As$ l'op\'erade des alg\`ebres associatives. M. Karoubi
a aussi d\'efini un foncteur \`a valeurs dans les alg\`ebres associatives (cf \cite {Ka1}, \cite {Ka2}). Mais ce foncteur
diff\`ere de notre foncteur $\Omega^{\mathcal As}$.
Enfin si on travaille avec $R$ un corps de caract\'eristique nulle on peut toujours d\'efinir un foncteur de formes
diff\'erentielles g\'en\'eralis\'ees pour $\ope$ une op\'erade unitaire, augment\'ee et born\'ee sup\'erieurement.
\end{subsection}

\begin{subsection}{Th\'eories de cocha\^{\i}nes et structures $\mathcal E_{\infty}$}
On fait d'abord quelques rappels sur les notions de th\'eories de cocha\^{\i}nes
\cite{Su}, \cite{Le} et de th\'eories cohomologiques
\cite{Ca}. 
On montre que tout foncteur de formes diff\'erentielles g\'en\'eralis\'ees est une th\'eorie de cocha\^{\i}nes. 
On explore 
la  structure de $\mathcal E_{\infty}$-alg\`ebre 
des th\'eories cohomologiques
 et des foncteurs de formes diff\'erentielles g\'en\'eralis\'ees; 
on en d\'eduit deux approches pour construire des cup-i produits sur
ces objets. Pour finir on montre que le foncteur des formes diff\'erentielles g\'en\'eralis\'ees pour les
$\mathcal E_{\infty}$-alg\`ebres est en un certain sens, universel.

On donne une version $\zed$-gradu\'ee des notions de th\'eories de cocha\^{\i}nes
et de th\'eories cohomologiques.

\begin{definition}
Soit $F:{\bf \mathcal S}^{op}\longrightarrow \mathcal O-\mathcal Alg_{dg}$ un foncteur. On dit que ce foncteur est une
th\'eorie de cocha\^{\i}nes s'il satisfait aux propri\'et\'es suivantes:

i)\qua Pour tout ensemble simplicial $X$ on a $\pi^*F(X)\cong H^*(X;R)$.

ii)\qua Pour toute inclusion simpliciale $i:K\longrightarrow X$, le morphisme induit 

$F(i):F(X)\longrightarrow F(K)$ est un \'epimorphisme.
\end{definition}

\begin{proposition}
Le foncteur des formes diff\'erentielles g\'en\'eralis\'ees ${\Omega^{\ope}}$ d\'efinit une th\'eorie de cocha\^{\i}nes.
\end{proposition}

\begin{proof}[Preuve]
La condition i) a \'et\'e v\'erifi\'ee dans le paragraphe pr\'ec\'edent (th\'eor\`eme 3.1).
V\'erifions que ${\Omega^{\ope}}$ transforme les
inclusions simpliciales en \'epimorphismes.

C'est une cons\'equence imm\'ediate du fait que ${\Omega^{\ope}}_*$ est un ensemble 
simplicial contractile.
En effet, on consid\`ere le diagramme d'ensembles simpliciaux suivants:
\begin{diagram}
K &             &\\
\dTo^{i}          & \rdTo^{f} &\\
X         & \rTo^{f'}  & {\Omega^{\ope}}_*\\  
\end{diagram}
Comme $i$ est une inclusion simpliciale (une cofibration dans la cat\'egorie des ensembles simpliciaux), 
et que  ${\Omega^{\ope}}_*$ 
est contractile et fibrant (comme ensemble simplicial), on peut toujours \'etendre une application simpliciale 
$f:K\longrightarrow {\Omega^{\ope}}_*$ en un morphisme
$f':X\longrightarrow {\Omega^{\ope}}_*$.
\end{proof}

Les axiomes dus \`a Cartan \cite{Ca} et \`a Swan \cite{Sw} permettent de construire une vaste classe de
th\'eories de cocha\^{\i}nes. Nous donnons ici une version $\zed$-gradu\'ee de ces axiomes.

\begin{definition}
Soit $M_*$ un $R$-module diff\'erentiel $\zed$-gradu\'e simplicial et augment\'e. On consid\`ere le foncteur
$$\mathcal M:{\bf \mathcal S}^{op}\longrightarrow \Mdg$$ 
tel que:
$$\mathcal M_X=Hom_{\bf \mathcal S}(X,M_*).$$ 
On dit que le foncteur $\mathcal M$ v\'erifie les axiomes de Cartan-Swan s'il satisfait les deux conditions suivantes:

i)\qua Le $R$-module $M_*$ est acyclique relativement \`a l'augmentation. Et le noyau de la diff\'erentielle
$d:M_*^0\longrightarrow M_*^1$ a le type d'homotopie d'un $K(R,0)$.

ii)\qua  Le $R$-module simplicial $M_*$ est un ensemble simplicial contractile.
\end{definition}

Le th\'eor\`eme suivant est d\^u \`a Cartan \cite{Ca} dans le cas $\mathbb N$-gradu\'e; on le
g\'en\'eralise au cas $\mathbb Z$-gradu\'e:

\begin{theorem}
Si le foncteur $\mathcal M$ associ\'e au $R$-module simplicial diff\'erentiel $\zed$-gradu\'e $M_*$
satisfait les axiomes de Cartan-Swan, alors il d\'efinit une th\'eorie de cocha\^{\i}nes. 
\end{theorem}

\begin{proof}[Preuve]
Cette g\'en\'eralisation ne pose pas de difficult\'es. Rappelons juste les arguments de Cartan.

Soit $Z^nM_*$ le noyau de la diff\'erentielle 
$$d:M^{n}_*\longrightarrow M^{n+1}_*$$ 
Comme le $R$-module $M_*$ est
acyclique relativement \`a l'augmentation, on a les suites exactes courtes:
$$0\longrightarrow Z^nM_*\longrightarrow M^{n}_*\stackrel{d}{\longrightarrow} Z^{n+1}M_*\longrightarrow 0$$
\noindent pour $n\geq 0$.

En tant que suites exactes courtes de groupes ab\'eliens simpliciaux, 
le morphisme $d$ est une fibration de Kan de fibre $Z^nM_*$. De plus comme $M^{n}_*$ est un ensemble 
simplicial contractile
et que $Z^0M_*$ a le type d'homotopie d'un $K(R,0)$, 
on en d\'eduit que chaque $Z^nM_*$ a le type d'homotopie d'un $K(R,n)$. Ces suites exactes s'identifient
aux fibrations:
$$K(R,n)\longrightarrow PK(R,n)\longrightarrow K(R,n+1).$$
Pour conclure on identifie $\pi^nM(X)$ avec $[X,Z^nM_*]$.
\end{proof}

Les $Z^nM_*$ forment un spectre dans la cat\'egorie des ensembles simpliciaux (en l'occurence un spectre
d'Eilenberg-Mac-Lane $HR$).
\begin{proposition}
Le foncteur des formes diff\'erentielles g\'en\'eralis\'ees ${\Omega^{\ope}}$ 
v\'erifie les axiomes de Cartan-Swan. On a les isomorphismes suivants:
$$\pi^n{\Omega^{\ope}}(X)\cong H^n(X,R)\cong [X,Z_n{\Omega^{\ope}}^*].$$
\end{proposition}

\begin{proof}[Preuve]
On v\'erifie que le noyau  $Z_0{\Omega^{\ope}}_*$ de la
diff\'erentielle
$$d:{\Omega^{\ope}}_*^0\longrightarrow {\Omega^{\ope}}_*^1$$ 
a le type d'homotopie d'un $K(R,0)$.

Dans la preuve du th\'eor\`eme pr\'ec\'edent $\pi_n{\Omega^{\ope}}_X$ est identifi\'e avec
$[X,Z_n{\Omega^{\ope}}_*]$, o\`u $Z^n{\Omega^{\ope}}_*$ est le noyau de
$$d:{\Omega^{\ope}}_*^n\longrightarrow {\Omega^{\ope}}_*^{n+1}.$$
De plus, on sait que le foncteur des formes diff\'erentielles g\'en\'eralis\'ees $\pi^n{\Omega^{\ope}}_X$
est isomorphe \`a $H^n(X;R)$. 

Ce qui nous permet l'identification de $\pi^0{\Omega^{\ope}}(X)$ avec [X,$Z_0{\Omega^{\ope}}_*$].
\end{proof}

On suppose maintenant que $\mathcal O=\mathcal E_{\infty}$ est un mod\`ele cofibrant de $\mathcal Com$ dans la
cat\'egorie des op\'erades unitaires augment\'ees.

Un des int\'er\^ets majeurs des $\mathcal E_{\infty}$-alg\`ebres est qu'elles apparaissent de mani\`ere naturelle
dans l'\'etude de l'homotopie des espaces topologiques.

En effet, Hinich et Schechtmann \cite{HS} ont d\'emontr\'e que pour tout ensemble simplicial $X$,
l'alg\`ebre des cocha\^{\i}nes singuli\`eres normalis\'ees est
munie de mani\`ere naturelle d'une structure de $\mathcal E_{\infty}$-alg\`ebre.

En coefficients $\mathbb F_p$, les op\'erations de Steenrod
sur $\pi^*C^*(X;\mathbb F_p)\cong H^*(X;\mathbb F_p)$ d\'etermin\'ees par cette structure $\mathcal E_{\infty}$
coincident avec les op\'erations de Steenrod classiques.

\begin{proposition}
Il existe un morphisme naturel de $\mathcal E_{\infty}$-alg\`ebres:
$${\Omega^{\mathcal E_{\infty}}}(X)\longrightarrow C^*(X;R)$$
\noindent qui induit une \'equivalence d'homotopie dans la cat\'egorie $\Mdg$.
\end{proposition}

\begin{proof}[Preuve]
On g\'en\'eralise la construction donn\'ee par M. Karoubi dans le cadre des formes diff\'erentielles
non-commutatives.

Pour ce faire on rappelle qu'une cocha\^{\i}ne singuli\`ere $\lambda\in C^n(\Delta[s];R)$ peut \^etre
consid\'er\'ee comme une application qui associe un \'el\'ement $\lambda(i_0,\ldots ,i_n)$ de $R$
\`a toute suite $(i_0,\ldots ,i_n)$ d'\'el\'ements de $\{0,\ldots,s\}$ (ce morphisme
doit aussi v\'erifier des conditions de compatibilit\'e avec les morphismes de faces et de
d\'eg\'en\'erescences de $\Delta[s]$).

Le cup produit 
$$C^n(\Delta[s];R)\times C^m(\Delta[s];R)\longrightarrow C^{n+m}(\Delta[s];R)$$
est donn\'e par la formule d'Alexander-Whitney:
$$(\lambda \cup \mu)(i_0,\ldots,i_{n+m})=\lambda(i_0,\ldots,i_n)\mu(i_n,\ldots,i_{n+m}).$$  
Si $\lambda\in C^0(\Delta[s];R)$, son bord $\delta(\lambda)$ est donn\'e par:
$\delta(\lambda)(i,j)=\lambda(i)-\lambda(j)$.

Nous allons d\'efinir un morphisme de $\mathcal E_{\infty}$-alg\`ebres:
$$\Phi_s:{\Omega^{\mathcal E_{\infty}}}(\Delta[s])\longrightarrow C^*(\Delta[s];R).$$
Comme 
${\Omega^{\mathcal E_{\infty}}}(\Delta[s])=\mathcal E_{\infty}(t_0,\ldots,t_s;dt_0,\ldots,dt_s)/I_s$, pour 
d\'efinir $\Phi_s$ il suffit de donner l'image de $\{t_0,\ldots,t_s\,dt_0,\ldots,dt_s\}$.

On pose $\Phi_s(t_r)=X_r$ avec $X_r\in R[X_0,\ldots,X_s]/(\sum_{r=0}^n X_r=1)$ ce polyn\^ome 
d\'efinit un \'el\'ement de  $C^0(\Delta[s];R)$: pour $i\in{0,\dots,s}$ $X_r(i)$ correspond 
\`a l'\'evaluation de ce polyn\^ome en 
(0,\ldots,0,1,0,\ldots,0) o\`u $1$ est en $i^{eme}$ position.

On d\'efinit $\delta X_r$ par $\delta X_r(i,j)=X_r(i)-X_r(j)$, enfin on
pose $\Phi_s(dt_r)=\delta X_r$.

Ces morphismes induisent des quasi-isomorphismes 
(en effet on remarque que l'image d'un scalaire $r\in R$ par $\Phi_s$ est l'application constante de valeur
$r$).

Tout ce qui pr\'ec\`ede permet de d\'efinir une transformation naturelle:
$$\Phi:{\Omega^{\mathcal E_{\infty}}}\longrightarrow C^*(-;R).$$
Et, par un argument de la th\'eorie des mod\`eles acycliques, on en d\'eduit que $\Phi$ induit
une \'equivalence d'homotopie dans la cat\'egorie des $R$-modules diff\'erentiels gradu\'es.
\end{proof}

On donne maintenant des r\'esultats de comparaison entre les diff\'erentes th\'eories de formes diff\'erentielles
g\'en\'eralis\'ees. Pour ce faire, on remarque que si $\ope$ est une op\'erade cofibrante unitaire augment\'ee,
alors toute $\mathcal E_{\infty}$-alg\`ebre est de
mani\`ere naturelle une $\ope$-alg\`ebre; cette  
structure est induite par le morphisme d'op\'erades $f$:
\begin{diagram}
 &             &    \mathcal E_{\infty}        \\
          & \ruTo^{f} &\dTo^{\epsilon}\\
\ope         & \rTo^{\epsilon}  & {\mathcal Com}\\
\end{diagram}
Le morphisme $f$ est un rel\`evement de l'augmentation $\epsilon$.
Celui-ci existe car $\ope$ est cofibrante et
l'augmentation $\epsilon:\mathcal E_{\infty}\longrightarrow \mathcal Com$ \'etant
une fibration triviale poss\`ede la propri\'et\'e de rel\`evement par rapport aux
cofibrations.

\begin{proposition}
Soit $\ope$ une op\'erade (cofibrante unitaire augment\'ee),
il existe une structure naturelle de $\ope$-alg\`ebre sur  $C^*(X;R)$ et
sur ${\Omega^{\mathcal E_{\infty}}}(X)$.

De plus, on a des morphismes naturels de $\ope$-alg\`ebres:
$${\Omega^{\ope}}(X)\longrightarrow{\Omega^{\mathcal E_{\infty}}}(X)$$
$${\Omega^{\ope}}(X)\longrightarrow C^*(X;R).$$
qui induisent des \'equivalences d'homotopie.
\end{proposition}

\begin{proof}[Preuve]
La structure de $\ope$-alg\`ebre sur $C^*(X;R)$ provient du morphisme 
$$f:\ope\longrightarrow \mathcal E_{\infty}$$ 
et du fait que $C^*(X;R)$ est une $\mathcal E_{\infty}$-alg\`ebre.

Le premier morphisme 
$$\Psi_X:{\Omega^{\ope}}(X)\longrightarrow{\Omega^{\mathcal E_{\infty}}}(X)$$
est obtenu via les applications: 
$$\Psi_s:{\Omega^{\ope}}(\Delta[s])\longrightarrow{\Omega^{\mathcal E_{\infty}}}(\Delta[s]).$$
On rappelle que:
$$\Omega^{\ope}(\Delta[s])=\ope (t_0,\ldots,t_s;dt_0,\ldots,dt_s)/I_s$$
et que:
$${\Omega^{\mathcal E_{\infty}}}(\Delta[s])=\mathcal E_{\infty}(t'_0,\ldots,t'_s;dt'_0,\ldots,dt'_s)/I_s.$$
Alors $\Psi_s$ est enti\`erement determin\'ee par les conditions suivantes $\Psi_s(t_i)=t'_i$ 
et $\Psi_s(dt_i)=dt'_i$.
\end{proof}

On donne maintenant deux mani\`eres de construire des cup-i produits pour le foncteur ${\Omega^{\ope}}$.

La premi\`ere m\'ethode applique la th\'eorie des mod\`eles acycliques.

Le groupe sym\'etrique $\Sigma_p$ agit sur ${\Omega^{\ope}}(X)^{\otimes p}$ via 
l'application d'\'echange $T$.

Soit $\alpha\in \Sigma_p$ une permutation cyclique d'ordre $p$. On introduit les op\'erations suivantes: 
$$\tau=1-\alpha$$
$$\sigma=1+\alpha+\ldots +\alpha^{p-1}.$$
Elles v\'erifient $\tau\sigma=\sigma\tau=0$. Consid\'erons le produit $\mu_0$ d\'efini par la
compos\'ee
$$\mu_0:{\Omega^{\ope}}(X)^{\otimes p}\stackrel{s_p}{\longrightarrow}
\ope(p)\otimes{\Omega^{\ope}}(X)^{\otimes p}\stackrel{\theta_p}{\longrightarrow} {\Omega^{\ope}}(X)$$
Avec $s_p$ une section de l'augmentation $\epsilon:\ope(p)\longrightarrow R$, et $\theta_p$ l'action
de $\ope$ sur ${\Omega^{\ope}}(X)$. 

Le morphisme $\mu_0\tau$ est nul en homotopie (le produit \'etant commutatif en homotopie, d'apr\`es un
argument de la th\'eorie des mod\`eles acycliques). Il
existe un homomorphisme $\mu_1$ de degr\'e -1, tel que:
$$\mu_0\tau=\mu_1d+d\mu_1.$$ 
Maintenant on consid\`ere le morphisme $\mu_1\sigma$: lui aussi est nul en homotopie. 
Toujours d'apr\`es la th\'eorie des mod\`eles acycliques il existe un homomorphisme $\mu_2$,
de degr\'e $-2$, tel que:
$$\mu_1\sigma=\mu_2d-d\mu_2.$$ 
Ainsi par it\'eration on construit une suite d'op\'erations v\'erifiant les formules suivantes:
$$\mu_{2n}\tau=\mu_{2n+1}d+d\mu_{2n+1}$$
$$\mu_{2n+1}\sigma=\mu_{2n}d-d\mu_{2n}.$$
Si on travaille \`a coefficients dans $\mathbb F_p$, les op\'erations de Steenrod sont obtenues 
en composant avec la diagonale:
$${\Omega^{\ope}}(X)\stackrel{\Delta}{\longrightarrow}{\Omega^{\ope}}(X)^{\otimes p}
\stackrel{\mu_n}{\longrightarrow}{\Omega^{\ope}}(X).$$ 
Bien s\^ur toutes ces constructions sont naturelles en $X$.

D\'ecrivons maintenant la seconde approche:

\begin{proposition}
Pour toute op\'erade cofibrante le foncteur ${\Omega^{\ope}}$ est \`a valeurs dans la cat\'egorie 
des $\mathcal E_{\infty}$-alg\`ebres.
\end{proposition}

\begin{proof}[Preuve]
On commence par montrer que ${\Omega^{\ope}}(X)$ est une alg\`ebre sur une op\'er\-ade acyclique que l'on note
$\mathcal End(\ope)$. 

On d\'efinit cette op\'erade par: 
$$\mathcal End(\ope)(n)=Hom_{\Mdg_{\bf \mathcal S}}({{\Omega^{\ope}}_*}^{\otimes n}, {\Omega^{\ope}}_*)$$ 
pour $n\geq 2$, on pose $\mathcal End(\ope)(0)=\mathcal End(\ope)(1)=R$.
Cette op\'erade est acyclique car on a prouv\'e que ${\Omega^{\ope}}_*$ est contractile comme R-module
simplicial diff\'erentiel $\zed$-gradu\'e.

Dans la cat\'egorie des op\'erades il existe un mod\`ele cofibrant not\'e $\mathcal E_{\infty}$ de
l'op\'er\-ade $\mathcal End(\ope)$:
$$\mathcal E_{\infty}\longrightarrow\mathcal End(\ope)\longrightarrow\mathcal Com.$$
Chaque morphisme est une fibration triviale donc $\mathcal E_{\infty}$ est aussi un mod\`ele 
cofibrant de l'op\'erade $\mathcal Com$.

Et le morphisme d'op\'erades $\mathcal E_{\infty}\longrightarrow\mathcal End(\ope)$ 
permet de d\'efinir une structure de $\mathcal E_{\infty}$-alg\`ebre sur ${\Omega^{\ope}}(X)$.
\end{proof}

\noindent {\bf Remarques}\qua Le $R$-module ${\Omega^{\ope}}(X)$ est aussi une alg\`ebre sur une seconde op\'erade not\'ee 
$TN_{\mathcal Dold}(\ope)$. Cette op\'erade est telle que:
$$TN_{\mathcal Dold}(\ope)(n)=Hom({\Omega^{\ope}}^{\otimes n},\Omega^{\ope})$$
pour $n\geq 2$, et
$$TN_{\mathcal Dold}(\ope)(0)=TN_{\mathcal Dold}(\ope)(1)=R.$$ 
Ce sont les transformations naturelles entre le foncteur
${\Omega^{\ope}}^{\otimes n}$ et le foncteur $\Omega^{\ope}$.  L'acyclicit\'e
de cette op\'erade est une cons\'equence imm\'ediate de la th\'eorie des
mod\`eles acycliques. Ceci g\'en\'eralise un r\'esultat de A. Dold \cite {Do}
obtenu pour le foncteur des cha\^{\i}nes singuli\`eres d'un ensemble
simplicial. De plus, il existe une paire de morphismes $\phi:TN_{\mathcal
Dold}(\ope)\rightleftharpoons \mathcal End(\ope):\psi$, telle que $\phi$  est
une fibration acyclique (dans la cat\'egorie des op\'erades) et $\phi\psi=Id$.

Soit $T$ un \'el\'ement de $TN_{\mathcal Dold}(\ope)(n)$, le morphisme
$\phi(T)\in \mathcal End(\ope)(n)$ est donn\'e par $\phi(T)=T(\Delta)$.On
rappelle que $\Delta$ est l'ensemble cosimplicial form\'e  par les
$\Delta[p]$.

Consid\'erons $f\in\mathcal End(\ope)(n)$, on lui associe la transformation naturelle $T_f$ d\'efinie comme
\'etant la compos\'ee:
$$T_f(X):Hom_{\mathbf S}(X,\Omega^{\ope}_*)^{\otimes n}\longrightarrow 
Hom_{\mathbf S}(X,{\Omega^{\ope}}^{\otimes n}_*)\stackrel{f^*}{\longrightarrow} 
Hom_{\mathbf S}(X,\Omega^{\ope}_*).$$
Les op\'erades $TN_{\mathcal Dold}(\ope)$ et $\mathcal End(\ope)$ ne sont pas cofibrantes car 
l'op\'erade $\mathcal Com$ est un r\'etract de ces deux op\'erades.

On r\'esume ces r\'esultats par le th\'eor\`eme suivant:

\begin{theorem}
Pour tout ensemble simplicial $X$ l'alg\`ebre $\Omega^{\mathcal O}(X)$ est munie de l'action
d'op\'erations. Cette action induit sur $\pi^*\Omega^{\mathcal O}(X)$ une structure d'alg\`ebre
instable sur l'alg\`ebre de Steenrod telle que  $\pi^*\Omega^{\mathcal O}(X)\cong
H^*(X,\mathbb F_p)$ soit
un isomorphisme d'alg\`ebres instables.
\end{theorem}

La proposition suivante relie les diff\'erents foncteurs de formes diff\'erentielles
g\'en\'eralis\'ees.

\begin{proposition}
{\rm i)}\qua Il existe un morphisme naturel de $\mathcal E_{\infty}$-alg\`ebres:
$${\Omega^{\mathcal E_{\infty}}}(X)\longrightarrow {\Omega^{\ope}}(X)$$
qui induit une \'equivalence d'homotopie dans la cat\'egorie $\Mdg$.

{\rm ii)}\qua Si $R$ est un corps de caract\'eristique nulle, on note $\mathcal A_{PL}$ le foncteur de Sullivan 
(\cite{Su}, \cite{Bo}, \cite{FHT1}). 
Alors on a  un morphisme naturel 
de $\mathcal E_{\infty}$-alg\`ebres:
$${\Omega^{\mathcal E_{\infty}}}(X)\longrightarrow \mathcal A_{PL}(X)$$
qui induit une \'equivalence d'homotopie dans la cat\'egorie $\Mdg$.

{\rm iii)}\qua Si on note $\Omega_{\mathcal K}$ le foncteur des formes diff\'erentielles non-commutatives  
introduit par M. Karoubi \cite{Ka1}, \cite{Ka2}, on a 
un morphisme naturel de $\mathcal E_{\infty}$-alg\`ebres:
$${\Omega^{\mathcal E_{\infty}}}(X)\longrightarrow {\Omega_{\mathcal K}}(X)$$
qui induit une \'equivalence d'homotopie dans la cat\'egorie $\Mdg$.
\end{proposition}

\begin{proof}[Preuve]
Comme pour les propositions pr\'ec\'edentes, on peut donner des formules explicites pour chacun de ces morphismes.  
\end{proof}
 
\end{subsection}

\begin{subsection}{Applications \`a l'homotopie des alg\`ebres sur une op\'erade}
On suppose que l'op\'erade $\ope$ permet de d\'efinir un foncteur de formes diff\'erent\-ielles
g\'en\'eralis\'ees. Dans ce cadre on va g\'en\'eraliser quelques constructions bien connues en homotopie
rationnelle. En fait on retrouve les constructions de l'homot\-opie rationnelle, quand on prend $R=\mathbb Q$
et $\mathcal O=\mathcal Com$, le foncteur $\Omega^{\mathcal Com}$ est le foncteur $A_{PL}$ de Sullivan
(\cite {Su},\cite {Bo}).

\begin{definition}
Pour deux $\ope$-alg\`ebres A et B on d\'efinit un espace fonctionnel simplicial $Hom^{\Delta}(A,B)$ dont les
 n-simplexes sont donn\'es par:
 
$Hom^{\Delta[n]}(A,B)=Hom_{\bf \ope-\mathcal Alg_{dg}}(A,\ope^n\amalg B)$.
\end{definition}

\begin{definition}
Le foncteur  
$$\mathbb M: {\ope-\mathcal Alg_{dg}^{op}} \longrightarrow \mathcal S$$
donn\'e par:
$$\mathbb M(A)=Hom_{\bf {\ope-\mathcal Alg_{dg}}}(A,{\Omega}^{\ope}_*)=Hom^{\Delta}(A,R).$$ 
est appel\'e foncteur de r\'ealisation simpliciale.
\end{definition}

\begin{proposition} Les foncteurs contravariants $\Omega^{\ope}$ et $\mathbb M$ sont des 
foncteurs adjoints au sens de Quillen 
entre les cat\'egories ${\bf \ope-\mathcal Alg_{dg}}$ et {\bf $\mathcal S$}.
\end{proposition}

\begin{proof}[Preuve]
On commence par montrer que $\Omega^{\ope}$ et $\mathbb M$ sont adjoints, c'est-\`a dire
que l'on a une bijection naturelle 
$$\Psi_X:Hom_{\bf \mathcal S}(X,Hom^{\Delta}(A,R))\longrightarrow 
Hom_{\bf {\ope-\mathcal Alg_{dg}}}(A,{\Omega}^{\ope}(X))$$
Si $X=\Delta[n]$ c'est \'evident (d'apr\`es la d\'efinition du foncteur $\mathbb M$).

Tout ensemble simplicial peut s'\'ecrire sous la forme d'un co\'egalisateur:
\begin{diagram}
\amalg_{s\in S}\Delta[s] &\pile{\rTo^f \\ \rTo_g} & \amalg_{t\in T}\Delta[t]& \rTo &X
\end{diagram}
On utilise la propri\'et\'e suivante du foncteur $\Omega^{\ope}$:

Le foncteur $\Omega^{\ope}:{\bf \mathcal S}\longrightarrow {\bf \ope-\mathcal Alg_{dg}}$ transforme les colimites en limites.
D'o\`u en appliquant $\Omega^{\ope}$ au co\'egalisateur pr\'ec\'edent on obtient l'\'egalisateur:
\begin{diagram}
\Omega^{\ope}(X) & \rTo  & \Omega^{\ope}(\amalg_{t\in T}\Delta[t]) &\pile{\rTo^{\Omega^{\ope}f} \\ 
\rTo_{\Omega^{\ope}g}} & \Omega^{\ope}(\amalg_{s\in S}\Delta[s])
\end{diagram}
Puis
\begin{diagram}
\Omega^{\ope}(X)& \rTo &\prod_{t\in T} \Omega^{\ope}(\Delta[t]) 
&\pile{\rTo \\ \rTo} & \prod_{s\in S} \Omega^{\ope}(\Delta[s])
\end{diagram}
et
\begin{diagram}
\Omega^{\ope}(X)\amalg R& \rTo &\prod_{t\in T} (\Omega^{\ope}(\Delta[t])\amalg R) 
&\pile{\rTo \\ \rTo} & \prod_{s\in S} (\Omega^{\ope}(\Delta[s]) \amalg R)
\end{diagram}
Le foncteur $Hom_{\bf \ope-\mathcal Alg_{dg}}(A,-)$ transforme les
produits en produits et les \'egal\-isat\-eurs en \'egalisateurs. On  
a alors:
\begin{diagram}
Hom_{\bf \ope-\mathcal Alg_{dg}}(A,\Omega^{\ope}(X)\amalg R)& \rTo 
&\prod_{t\in T}Hom_{\bf \ope-\mathcal Alg_{dg}}(A, \Omega^{\ope}(\Delta[t])\amalg R)\\ 
\ldots &\pile{\rTo \\ \rTo} & \prod_{s\in S} Hom_{\bf \ope_-\mathcal Alg{dg}}(A,\Omega^{\ope}(\Delta[s]) \amalg R)
\end{diagram}
Si on applique le foncteur $Hom_{\bf \mathcal S}(-,Hom^{\Delta}(A,B))$ au co\'egalisateur initial, 
on obtient l'\'egalisateur:
\begin{diagram}
Hom_{\bf \mathcal S}(X,Hom^{\Delta}(A,R)) &\rTo& \prod_{t\in T}Hom_{\bf \mathcal S}(\Delta[t],Hom^{\Delta}(A,R))\\
\ldots &\pile{\rTo \\ \rTo} & \prod_{s\in S}Hom_{\bf \mathcal S}(\Delta[s],Hom^{\Delta}(A,R))
\end{diagram}
Comme $\Psi$ est une bijection sur les seconds et troisi\`emes termes, c'est aussi une bijection sur les premiers. 

Pour finir, on remarque que le foncteur $\Omega^{\ope}$ transforme les cofibrations entre
ensembles simpliciaux en fibrations et qu'il pr\'eserve les \'equivalences faibles.
\end{proof}

Gr\^ace \`a la notion de formes diff\'erentielles g\'en\'eralis\'ees, on peut d\'efinir un objet chemin
pour la cat\'egorie des $\ope$-alg\`ebres: en effet consid\'erons
la $\ope$-alg\`ebre $\Omega^{\ope}_1=\Omega^{\ope} (\Delta [1])$ (avec $\Delta[1]$ le 1-simplexe standard). 
Puis on d\'efinit le foncteur:
$$-\amalg\Omega^{\ope}_1:{\bf \ope-\mathcal Alg_{dg}}\longrightarrow {\bf \ope-\mathcal Alg_{dg}}$$
Pour toute $\ope$-alg\`ebre $A$ on a les applications suivantes:
$$A\stackrel{p}{\longrightarrow}A\amalg\Omega^{\ope}_1\stackrel{(d_0,d_1)}{\longrightarrow}  A\times A.$$
On rappelle que l'on peut identifier $\Omega^{\ope}_1$ avec $\ope(t,dt)$,
on a des morphismes
$$\delta_0,\delta_1:\Omega^{\ope}_1\longrightarrow \Omega^{\ope}_0\cong R$$
tels que $\delta_i(dt)=0$ et $\delta_i(t)=i$.
Le morphisme p est l'application canonique.
Et les applications $d_0$, $d_1$ sont d\'efinies comme \'etant l'identit\'e sur $A$ et $\delta_i$ 
sur $\Omega^{\ope}_1$.

Ces morphismes satisfont les propri\'et\'es suivantes:

i)\qua $d_ip$ est l'identit\'e, $i=0,1$.

ii)\qua $d_i$ est une fibration triviale, $i=0,1$. Le morphsime $(d_0,d_1)$ est une fibration.

iii)\qua $p$ est une \'equivalence faible.

Gr\^ace \`a ces observations on d\'efinit maintenant une notion d'homotopie simpliciale.

\begin{definition}
Deux morphismes $f,g:A\longrightarrow B$ sont simplicialement homotopes s'il existe une application
H de $Hom^{\Delta[1]}(A,B)$, $H:A \longrightarrow B\amalg \Omega^{\ope}_1$ telle que $d_0H=f$ et $d_1H=g$. 
\end{definition}

\begin{proposition}
L'homotopie simpliciale satisfait les propri\'et\'es suivantes:

{\rm i)}\qua Soient deux applications $f,g:A\longrightarrow B$ qui sont simplicialement homotopes, alors:   
 $$\pi^*f=\pi^*g:\pi^*A\longrightarrow \pi^*B.$$
{\rm ii)}\qua Si $A$ est une $\ope$-alg\`ebre cofibrante l'homotopie simpliciale
est une relation d'\'equivalence.

{\rm iii)}\qua Si $A$ est cofibrante on a l'isomorphisme:
$$[A,B]_{Ho{\bf \ope-\mathcal Alg_{dg}}}\cong\pi^0Hom^{\Delta}(A,B).$$
\end{proposition}

\end{subsection}

\end{section}

\begin{section}{Un mod\`ele de la suite spectrale de Leray-Serre}

Un des outils fondamentaux en homotopie rationnelle est la construction du
mod\`ele alg\'ebrique de la fibre \`a partir de celui de la fibration. Construction
que l'on doit \`a P.P. Grivel \cite {G} (pour les espaces 1-connexes) et \`a S. Halperin \cite {Ha}
(pour les espaces nilpotents).

Nous g\'en\'eralisons l'approche de P.P. Grivel pour le foncteur $\mathcal A_{PL}$ de Sullivan
aux formes diff\'erentielles g\'en\'eralis\'ees pour les $\mathcal O$-alg\`ebres.

Les quatre premiers paragraphes sont consacr\'es \`a la construction d'une suite spectrale de 
Leray-Serre au moyen des formes diff\'erentielles g\'en\'eralis\'ees. A cet effet, on d\'efinit
des formes diff\'erentielles g\'en\'eralis\'ees avec coefficients locaux (paragraphe $1$), puis
des formes diff\'erentielles g\'en\'eralis\'ees pour les bisimplexes (paragraphe $2$).

Les techniques bisimpliciales (Dress \cite {Dr},\cite {MC}) sont utilis\'ees pour associer un bicomplexe \`a une 
fibration (paragraphe $3$). En filtrant ce bicomplexe,
 on retrouve la suite spectrale de Leray-Serre (paragraphe $4.3$,
th\'eor\`eme $4.1$).

Dans le paragraphe $4.4$, on construit une suite spectrale \`a partir du mod\`ele cofibrant en
$\mathcal E_{\infty}$-alg\`ebres d'une fibration. Par un th\'eor\`eme de comparaison entre
cette nouvelle suite spectrale et la suite spectrale de Leray-Serre construite pr\'ec\'edemment,
on en d\'eduit un mod\`ele de la fibre (proposition $4.9$).

Le paragraphe $4.5$ donne une interpr\'etation alg\'ebrique (sur le mod\`ele) de la transgression
pour la suite spectrale de Leray-Serre.   

\begin{subsection}{Pr\'efaisceaux sur les ensembles simpliciaux}Soit $X$ un
ensemble simplicial; on lui associe une cat\'egorie que l'on note $\bf X$. Les
objets de $\bf X$ sont les $p$-simplexes de $X$. Les morphismes de $\bf X$
correspondent aux applications $\alpha:\Delta[p]\longrightarrow \Delta[q]$
telles que le diagramme ci-dessous commute: \begin{diagram}
\Delta[p]&&\rTo^{\alpha} &&\Delta[q]\\
&\rdTo^{\alpha^*x}    &              &\ldTo^{x}    &\\
         &&     X        &&         \\
\end{diagram}
A toute application simpliciale $f:X\longrightarrow Y$ est associ\'e un foncteur ${\bf f}:\bf X\longrightarrow \bf Y$.

Un pr\'efaisceau sur $X$ \`a valeurs dans une cat\'egorie $\mathcal C$ est un foncteur contravariant de $\bf X$ 
dans $\mathcal C$.
On note $\mathcal F(\alpha,x):\mathcal F(\alpha^*x)\longrightarrow\mathcal F(x)$ l'application induite
par $\alpha:\Delta[p]\longrightarrow \Delta[q]$.

Les pr\'efaisceaux que l'on consid\`ere sont \`a valeurs dans $\bf Ab$ (la cat\'egorie des groupes ab\'eliens),
dans $\Mdg$ ou dans $\mathcal O-\mathcal Alg_{dg}$.

Si $f:X\longrightarrow Y$ est une application simpliciale et si $\mathcal F$ est un pr\'efaisceau sur $X$ alors on
a un pr\'efaisceau image r\'eciproque $f^*\mathcal F=\mathcal F{\bf f}$.

Les pr\'efaisceaux sur $X$ \`a valeurs dans une cat\'egorie ab\'elienne forment encore une cat\'egorie ab\'elienne.
Une suite de pr\'efaisceaux 
$$0\longrightarrow \mathcal F\longrightarrow \mathcal G\longrightarrow \mathcal H \longrightarrow 0$$
est exacte si pour tout $p$-simplexe $x_p$ de $X$ la suite 
$$0\longrightarrow \mathcal F(x_p)\longrightarrow \mathcal G(x_p)\longrightarrow \mathcal H(x_p)\longrightarrow 0$$   
est exacte.

Les foncteurs formes diff\'erentielles g\'en\'eralis\'ees sont des pr\'efaisceaux.
Les co\-cha\^{\i}nes singuli\`eres sont aussi un exemple de pr\'efaisceau.

{\noindent \bf Pr\'efaisceaux simpliciaux}\qua Un exemple important est le cas des pr\'efaisceaux \`a valeurs dans
$\bf SC$ la cat\'egorie des objets
simpliciaux de $\bf C$. Soit $\mathcal F$ un pr\'efaisceau \`a valeurs dans $\bf SC$. 

Une section $s$ du 
pr\'efaisceau $\mathcal F$ est la donn\'ee pour tout $q$-simplexe $x\in X_q$ d'un ensemble simplicial $s(x)$ tel que
pour toute application $\alpha:[p]\longrightarrow [q]$ on a: 
$$\mathcal F(\alpha,x)(s(\alpha^*(x)))=\alpha^*(s(x)).$$
Un pr\'efaisceau constant $\mathcal F$ sur $X$ est tel que pour tout $q$-simplexe $x\in X_q$, $\mathcal F(x)=M$
avec $M$ un objet de  $\bf SC$; une section de $\mathcal F$ est une application simpliciale de $X$ dans $M$.

{\noindent \bf Syst\`emes de coefficents locaux}\qua Un syst\`eme de coefficients locaux $\mathcal L$ sur un ensemble simplicial $X$
est
un pr\'efaisceau sur $X$ \`a valeurs dans $\bf Ab$ tel que pour tout op\'erateur de face 
$\delta_i:\Delta[p+1]\longrightarrow \Delta[p]$ et pour tout $p$-simplexe $x\in X_p$, les morphismes
$\mathcal L(\delta_i;x):\mathcal L(\delta_ix)\longrightarrow\mathcal L(x)$ soient des isomorphismes.

Pour tout 0-simplexe $x\in X$, on a une action de $\pi_1(X,x)$ sur le groupe ab\'elien $\mathcal L(x)$.
Le syst\`eme de coefficients est dit simple si l'action des groupes fondamentaux est triviale. De plus
si $X$ est connexe et si le syst\`eme de coefficients est simple alors les groupes $\mathcal L(x)$ sont
isomorphes.

On construit un foncteur $\Omega^{\mathcal O}(X;-):\bf{F}\Mdg_X\longrightarrow \bf{SF}\Mdg_X$
de la cat\'egorie des pr\'efaisceaux sur $X$ \`a valeurs dans $\Mdg$ vers la cat\'egorie des pr\'efaisceaux sur
$X$ \`a valeurs dans ${\bf S}\Mdg$.

On fixe $\mathcal F$ un pr\'efaisceau sur $X$ \`a valeurs dans $\Mdg$. On lui associe
un pr\'efaisceau not\'e $\Omega^{\mathcal O}(X;\mathcal F)$ sur $X$ \`a valeurs dans $\bf S\Mdg$. Ce
pr\'efaisceau associe \`a $x\in X$ le module diff\'erentiel gradu\'e simplicial 
$\Omega^{\mathcal O}_*\otimes\mathcal F(x)$.

\begin{definition}
Une section du pr\'efaisceau sur $\Omega^{\mathcal O}(X;\mathcal F)$ s'appelle une forme 
diff\'erentielle g\'en\'eralis\'ee
sur $X$ \`a valeurs dans le pr\'efaisceau $\mathcal F$. On note\nl 
$\underline\Omega^{\mathcal O}(X;\mathcal F)$ l'ensemble de ses sections.
\end{definition}

L'objet $\underline\Omega^{\mathcal O}(X;\mathcal F)$ est naturellement bigradu\'e. 
L'ensemble des \'el\'ements de bi\-degr\'e $(r,s)$ est not\'e
$\underline\Omega^{\mathcal O,r}(X;\mathcal F^s)$. Le degr\'e $s$ est donn\'e par le degr\'e 
diff\'erentiel de $\mathcal F$, et le degr\'e $r$ par celui de $\Omega^{\mathcal O}_*$.
On a donc bien un foncteur
$\underline\Omega^{\mathcal O,*}(X;\mathcal F^*)$ de la cat\'egorie des ensembles simpliciaux
vers la cat\'egorie des $R$-modules diff\'erentiels bigradu\'es.

Int\'eressons-nous au cas o\`u $\mathcal F$ est un pr\'efaisceau constant. Alors
$\mathcal F$ s'identifie \`a un R-module diff\'erentiel gradu\'e $M$ concentr\'e en  degr\'e 0.
Une forme diff\'erent\-ielle
g\'en\'eralis\'ee \`a valeurs dans $\mathcal F$ est une application simpliciale de $X$ dans 
$\Omega^{\mathcal O}_*\otimes M$. En particulier si $M=R$ alors 
on a $\underline\Omega^{\mathcal O}(X;R)=\Omega^{\mathcal O}(X)$.

\begin{proposition} Soit $X$ un ensemble simplicial et $\mathcal F=M$ un pr\'efaisceau constant, avec $M$ un $R$-module
plat; alors on a l'isomorphisme:
$$\pi^*\underline\Omega^{\mathcal O}(X;M)\cong H^*(X;M).$$
\end{proposition}
\begin{proof}[Preuve]
On utilise la th\'eorie des mod\`eles acycliques.
Le foncteur $\underline\Omega^{\mathcal O}(-;M)$ est corepr\'esentable car l'ensemble simplicial
$\Omega^{\mathcal O}_*\otimes M$ est contractile.

Le foncteur est acyclique car pour tout $p$ le $R$-module diff\'erentiel gradu\'e 
$\Omega^{\mathcal O}(\Delta[p])\otimes M$
est acyclique. 
\end{proof}

De plus le foncteur $\underline\Omega^{\mathcal O}(X;-)$ est exact:

\begin{proposition} 
Si $0\longrightarrow \mathcal F\longrightarrow \mathcal G\longrightarrow \mathcal H\longrightarrow 0$ est une suite
exacte de pr\'efaisceaux, alors la suite 
$$0\longrightarrow \Omega^{\mathcal O}(X;\mathcal F)\longrightarrow 
\Omega^{\mathcal O}(X;\mathcal G)\longrightarrow\Omega^{\mathcal O}(X; \mathcal H)\longrightarrow 0$$
est exacte, tout comme la suite
$$0\longrightarrow \underline\Omega^{\mathcal O}(X;\mathcal F)\longrightarrow 
\underline\Omega^{\mathcal O}(X;\mathcal G)\longrightarrow
\underline\Omega^{\mathcal O}(X; \mathcal H)\longrightarrow 0.$$
\end{proposition}
\begin{proof}[Preuve]
On suppose que 
$0\longrightarrow \mathcal F\stackrel{u}{\longrightarrow} \mathcal G
\stackrel{v}{\longrightarrow} \mathcal H\longrightarrow 0$ 
est une suite exacte de pr\'efaisceaux sur un ensemble simplicial $X$ ce qui signifie que pour tout 
$x\in X$ la suite
$$0\longrightarrow \mathcal F(x)\stackrel {u_x}{\longrightarrow} 
\mathcal G(x)\stackrel {v_x}{\longrightarrow}\mathcal H(x)\longrightarrow 0$$
est exacte. La platitude de $\Omega(\Delta)$ entra\^{\i}ne que
$$0\longrightarrow \Omega^{\mathcal O}_*\otimes \mathcal F(x)\stackrel {Id\otimes u_x}{\longrightarrow} 
\Omega^{\mathcal O}_*\otimes\mathcal G(x)\stackrel {Id\otimes v_x}{\longrightarrow}
\Omega^{\mathcal O}_*\otimes\mathcal H(x)\longrightarrow 0$$
est une suite exacte courte. Donc 
$$0\longrightarrow \Omega^{\mathcal O}_*\otimes \mathcal F\stackrel {Id\otimes u}{\longrightarrow} 
\Omega^{\mathcal O}_*\otimes\mathcal G\stackrel {Id\otimes v}{\longrightarrow}
\Omega^{\mathcal O}_*\otimes\mathcal H\longrightarrow 0$$
est aussi une suite exacte courte.

L'exactitude de la suite 
$$0\longrightarrow \underline\Omega^{\mathcal O}(X;\mathcal F)\longrightarrow 
\underline\Omega^{\mathcal O}(X;\mathcal G)\longrightarrow
\underline\Omega^{\mathcal O}(X; \mathcal H)\longrightarrow 0.$$
se d\'eduit de la propri\'et\'e d'extension des formes diff\'erentielles.
\end{proof}

Soit $\mathcal F$ est un pr\'efaisceau sur $X$ \`a valeurs dans $\Mdg$.
On note $d$ sa diff\'erent\-ielle, $\mathcal Z^s$ le pr\'efaisceau des cobords 
de degr\'e $s$ de $\mathcal F$ et $\mathcal H^s(\mathcal F)$ le pr\'efaisceau des classes de cohomologie de
degr\'e $s$ de $\mathcal F$. De la suite exacte:
$$0\longrightarrow\mathcal Z^s\longrightarrow\mathcal F^s\longrightarrow \mathcal H^s(\mathcal
F)\longrightarrow 0$$ 
et de la proposition pr\'ec\'edente on d\'eduit le r\'esultat suivant:
 
\begin{coroll}
On a $\pi^s(\underline\Omega^{\mathcal O}(X;\mathcal F))\cong \underline
\Omega^{\mathcal O}(X;\mathcal H^s(\mathcal F))$ o\`u 
$\pi^s(\underline \Omega^{\mathcal O}(X;\mathcal F))$ est la cohomologie de
$\underline \Omega^{\mathcal O}(X;\mathcal F)$ calcul\'ee avec la diff\'erentielle de $\mathcal F$. 
\end{coroll}

\end{subsection}

\begin{subsection}{Les formes diff\'erentielles g\'en\'eralis\'ees pour les
bisimplexes}

Le $R$-module diff\'erentiel $\Omega^{\mathcal O}_*\otimes\Omega^{\mathcal O}_*$ est un $R$-module 
diff\'erentiel 
bisimplicial et bigradu\'e, un \'el\'ement de bidimension $(p,q)$ est un \'el\'ement de 
$\Omega^{\mathcal O}(\Delta[p])\otimes\Omega^{\mathcal O}(\Delta[q])$. On note $\Omega_{bi}^{\mathcal O}(X)$ 
l'ensemble des formes diff\'erentielles 
g\'en\'eralis\'ees pour le bisimplexe $X$, 
$\Omega_{bi}^{\mathcal O}(X)=Hom_{Bi\mathcal S}(X,\Omega^{\mathcal O}_*\otimes\Omega^{\mathcal O}_*)$.

Une forme diff\'erentielle g\'en\'eralis\'ee de bidimension $(p,q)$ sur un ensemble bisimplicial $X$
est une application bisimpliciale de $X$ dans $\Omega^{\mathcal O}_p\otimes\Omega^{\mathcal O}_q$.  

L'objet $\Omega_*^{\mathcal O}\otimes\Omega^{\mathcal O}_*$ poss\`ede une autre bigraduation donn\'ee par le
degr\'e diff\'erentiel.

Une forme diff\'erentielle g\'en\'eralis\'ee de bidegr\'e $(l,m)$ sur $X$ un ensemble bisimplicial
est une application bisimpliciale de $X$ dans $\Omega^{\mathcal O,l}_*\otimes\Omega^{\mathcal O,m}_*$.

On a un foncteur $\beta:\mathcal S\longrightarrow Bi\mathcal S$
de la cat\'egorie des ensembles 
simpliciaux dans la cat\'egorie des ensembles bisimpliciaux tel que:
$$(\beta X)_{p,q}=X_p$$
 pour tout $q$. On v\'erifie ais\'ement que pour tout ensemble simplicial $X$ on a:
$$\Omega_{bi}^{\mathcal O}(\beta X)=\Omega^{\mathcal O}(X).$$
Soit $f:E\longrightarrow B$ une surjection simpliciale. A ce morphisme on associe un ensemble bisimplicial
$S_{.,.}f$. Un \'el\'ement de $S_{p,q}f$ est un couple d'applications 

$w:\Delta[p]\times\Delta[q]\longrightarrow E$ et $u:\Delta[p]\longrightarrow B$ telles que:
\begin{diagram}
\Delta[p]\times\Delta[q] & \rTo^{w} & E \\
\dTo^{pr_1}& &\dTo^{f}\\
\Delta[p]& \rTo^{u} & B \\
\end{diagram}
On peut montrer que $S_{.,.}f$ et $\beta E$ sont quasi-isomorphes \cite {G}. 
Comme $\Omega_{bi}^{\mathcal O}(\beta E)$ $=\Omega^{\mathcal O}(E)$ on en d\'eduit  la proposition:

\begin{proposition} Soit $f:E\longrightarrow B$ une surjection simpliciale. On a l'isomor\-phisme:
$$\pi^*(\Omega_{bi}^{\mathcal O}(S_{.,.}f))\cong H^*(E;R)$$
\end{proposition}

Le bisimplexe $S_{.,.}f$ permet de d\'efinir un syst\`eme de coefficients locaux sur $B$.
On consid\`ere $B$ comme une cat\'egorie (les objets de $B$ sont les p-simplexes) et on
d\'efinit un foncteur $\mathbb F:B\longrightarrow \mathcal S$ en associant \`a tout
simplexe $u_p:\Delta[p]\longrightarrow B$ l'ensemble simplicial:
$$(\mathbb F(u_p))_q=\{w_{p,q}/(w_{p,q},u_p)\in S_{p,q}f\}$$
On a alors le syst\`eme de coefficients locaux $\mathcal F:B\longrightarrow \Mdg$ en posant
$\mathcal F(u_p)=\Omega^{\mathcal O}(\mathbb F(u_p))$. On a  un pr\'efaisceau
$\Omega^{\mathcal O}(B,\mathcal F)$ sur $B$ \`a valeurs dans $\mathcal F$. 

\begin{proposition}
On a $H^*(\underline\Omega^{\mathcal O,r}(B;\mathcal F^s);d_{\mathcal F})
\cong \underline\Omega^{\mathcal O,r}(B;\mathcal H^s(\mathcal F))$.
De plus si $f:E\longrightarrow B$ est une fibration de Kan et si le syst\`eme
de coefficients locaux $\mathcal F$ est simple, alors on a aussi:
$$H^*(\underline\Omega^{\mathcal O,r}(B;\mathcal F^s);d_{\mathcal F})\cong 
\underline\Omega^{\mathcal O,r}(B;H^s(\mathcal F)).$$
\end{proposition}

\end{subsection}

\begin{subsection}{Formes diff\'erentielles et suites spectrales de
Leray-Serre}

Dans ce paragraphe on donne une construction de la suite spectrale de Leray-Serre
\`a partir des formes diff\'erentielles g\'en\'eralis\'ees.

On utilise cette suite spectrale dans les paragraphes suivants dans le cadre des alg\`ebres sur une 
$E_{\infty}$-op\'erade cofibrante. Ceci permet de donner un mod\`ele alg\'ebrique des fibrations.

\begin{theorem}
Soit $f:E\longrightarrow B$ une application surjective entre ensembles simpliciaux (on suppose que la 
cohomologie de $B$ ou que $\mathcal H^s(\mathcal F)$ sont finis).

{\rm i)}\qua On a une suite spectrale qui converge vers $H^*(E;R)$, telle que:

$E^{r,s}_0=\underline\Omega^{\mathcal O,r}(B;{\mathcal F^s})$,
$E^{r,s}_1=\underline\Omega^{\mathcal O,r}(B;\mathcal H^s(\mathcal F))$ et
$E^{r,s}_2=H^r(B;\mathcal H^s(\mathcal F))$

o\`u $\mathcal F$ est le syst\`eme de coefficients locaux sur $B$ d\'efini pr\'ec\'edemment.

{\rm ii)}\qua On suppose que $f$ est une fibration de Kan, que $B$ est connexe et que le syst\`eme de coefficients 
locaux $\mathcal H^s(\mathcal F)$ est simple (ce qui est le cas si $B$ est $1$-connexe). Alors
$$E^{r,s}_2=H^r(B;H^s(F;R))$$
$F$ d\'esignant la fibre de $f$.
\end{theorem}
\proof[Preuve]
On consid\`ere la premi\`ere filtration du bicomplexe $\Omega_{bi}^{\mathcal O}(S_{r,s}f)$. Une forme 
$\omega\in \Omega_{bi}^{\mathcal O}(S_{r,s}f)$
est de filtration $n$, si pour tout bisimplexe $(w,u)\in S_{p,q}f$ on a 
$\omega(w,u)\in\bigoplus_{r\geq n}\Omega^{\mathcal O,r}_p\otimes\Omega^{\mathcal O}_q$.

De cette filtration on d\'eduit une suite spectrale telle que: 
$E^{r,s}_0=\underline\Omega^{\mathcal O,r}(B;{\mathcal F^s})$.
En fait on explicite deux morphismes de modules diff\'erentiels gradu\'es:
$$\phi:\underline\Omega^{\mathcal O,r}(B;{\mathcal F^s})\rightleftharpoons\Omega_{bi}^{\mathcal O}(S_{r,s}f):\psi$$
inverses l'un de l'autre.

L'application $\phi$ est d\'efinie de la mani\`ere suivante.

Soit $\omega\in \underline\Omega^{\mathcal O,r}(B;{\mathcal F^s})$ et $u_p$ un $p$-simplexe de $B$ on a
$$\omega(u_p)\in\Omega^{\mathcal O,r}_p\otimes\mathcal F^s(u_p)$$
et par d\'efinition
$$\Omega^{\mathcal O,r}_p\otimes\mathcal F^s(u_p)=\Omega^{\mathcal O,r}_p\otimes
\Omega^{\mathcal O,s}(\mathbb F(u_p))=\underline\Omega^{\mathcal O,s}(\mathbb F(u_p);\Omega^{\mathcal O,r}_p).$$
L'application bisimpliciale $\phi(\omega)$ envoie $(w_{p,q},u_p)\in S_{p,q}f$ sur $\omega(u_p)(w_{p,q}).$

A l'inverse, soit $\omega\in\Omega_{bi}^{\mathcal O}(S_{r,s}f)$ et $u_p$ un $p$-simplexe de $B$, alors
pour $w_{p,q}\in (\mathbb F(u_p))_q$ on pose 
$$\psi(\omega)(u_p)(w_{p,q})=\omega(w_{p,q},u_p)$$
Comme la diff\'erentielle de $E^{r,s}_0$ est induite par la deuxi\`eme diff\'erentielle de 
$\Omega_{bi}^{\mathcal O}(S_{r,s}f)$, elle s'identifie \`a la diff\'erentielle $d_{\mathcal F}$ 
sur $\underline\Omega^{\mathcal O,r}(B;{\mathcal F^s})$. On a donc l'isomorphisme:
$$H^*(\underline\Omega^{\mathcal O,r}(B;\mathcal F^s);Id\otimes d)
\cong \underline\Omega^{\mathcal O,r}(B;\mathcal H^s(\mathcal F))$$
d'o\`u:
$$E^{r,s}_1=\underline\Omega^{\mathcal O,r}(B;\mathcal H^s(\mathcal F))$$
et par suite:
$$E^{r,s}_2=H^r(B;\mathcal H^s(\mathcal F)).\eqno{\qed}$$
\end{subsection}

\begin{subsection}{Mod\`ele d'une fibration}

On suppose maintenant que l'op\'erade $\mathcal O=\mathcal E_{\infty}$ est un
mod\`ele cofibrant de l'op\'erade $\mathcal Com$.
L'existence de la suite spectrale de Leray-Serre associ\'ee aux formes diff\'erentielles
g\'en\'eralis\'ees pour les $\mathcal E_{\infty}$-alg\`ebres nous
permet de d\'ecrire un mod\`ele alg\'ebrique des fibrations. Au pr\'ealable nous avons besoin de
r\'esultats concernant le coproduit des $\mathcal E_{\infty}$-alg\`ebres.

{\bf \noindent Le produit tensoriel des op\'erades}\qua Soient $\mathcal P$ et $\mathcal Q$ deux op\'erades. On d\'efinit
l'op\'erade $\mathcal P\otimes \mathcal Q$ en posant $(\mathcal P\otimes Q)(n)=\mathcal P(n)\otimes Q(n)$. On remarque
que si $A$ est une $P$-alg\`ebre et si $B$ est une $Q$-alg\`ebre alors $A\otimes B$ est une $P\otimes Q$-alg\`ebre.

{\bf \noindent Le coproduit des $\mathcal E_{\infty}$-alg\`ebres}

\begin{proposition}
Il existe un morphisme d'op\'erades $\Delta:\mathcal E_{\infty}\longrightarrow \mathcal E_{\infty}\otimes \mathcal E_{\infty}$.
\end{proposition}

\begin{proof}[Preuve]
Consid\'erons le morphisme canonique $\Delta^c:\mathcal Com\longrightarrow \mathcal Com\otimes\mathcal Com$. Il existe
une application $\Delta$ telle que le diagramme ci-dessous commute:
\begin{diagram}
\mathcal E_{\infty} & \rDotsto^{\Delta}  & \mathcal E_{\infty}\otimes
\mathcal E_{\infty} \\
\dTo& &\dTo\\
\mathcal Com & \rTo^{\Delta^c}  & \mathcal Com\otimes
\mathcal Com.\\
\end{diagram}
L'existence de $\Delta$ est une cons\'equence de la propri\'et\'e de rel\`evement \`a gauche par rapport
\`a la fibration triviale
$\mathcal E_{\infty}\otimes\mathcal E{\infty}\longrightarrow\mathcal Com\otimes\mathcal Com$.
\end{proof}

\begin{coroll}Soient $A$ et $B$ deux $\mathcal E_{\infty}$-alg\`ebres alors $A\otimes B$ poss\`ede
une structure naturelle de $\mathcal E_{\infty}$-alg\`ebre.
\end{coroll}
\begin{proof}[Preuve]
L'objet $A\otimes B$ est une $\mathcal E_{\infty}\otimes\mathcal E_{\infty}$-alg\`ebre. En cons\'equence,
si on a un morphisme $\Delta:\mathcal E_{\infty}\longrightarrow\mathcal E_{\infty}\otimes\mathcal E_{\infty}$, alors
on peut utiliser la restriction de structure.
\end{proof}

On a des morphismes canoniques de $\mathcal E_{\infty}$-alg\`ebres
$A\longrightarrow A\otimes B$, $B\longrightarrow A\otimes B$. De la propri\'et\'e
universelle du coproduit, on obtient un morphisme:
$$p:A\amalg B\longrightarrow A\otimes B.$$
On a le r\'esultat suivant:

\begin{proposition}
Soient $A$ et $B$ deux $\mathcal E_{\infty}$-alg\`ebres cofibrantes,
alors le morphisme ci-dessus:
$$p:A\amalg B\longrightarrow A\otimes B$$
est une fibration triviale.
\end{proposition}

\begin{proof}[Preuve]
Pour commencer, on consid\`ere le morphisme de $R$-modules diff\'erent\-iels gradu\'es:
$$s:A\otimes B\longrightarrow A\amalg B$$
tel que $s(a\otimes b)=1_2\otimes a\otimes b$. V. Hinich a  d\'emontr\'e \cite {Hi2} que $s$ est un quasi-isomorphisme.
On v\'erifie ais\'ement que:
$p(s(a\otimes b))=a\otimes b$. La propri\'et\'e s'ensuit.
\end{proof}

La proposition suivante est une cons\'equence directe de la proposition $4.5$:

\begin{proposition}
Le foncteur $\Omega^{\mathcal E_{\infty}}_{bi}$ est un foncteur de la cat\'egorie
des ensembles bisimpliciaux dans la cat\'egorie des $\mathcal E_{\infty}$-alg\`ebres.
\end{proposition}

On a \'egalement:

\begin{proposition}
Soient $A$ une $\mathcal E_{\infty}$-alg\`ebre
quasi-libre 1-connexe et $A\longrightarrow A\amalg_{\tau} \mathcal E_{\infty}(V)$ un morphisme quasi-libre. Alors
on a une  suite spectrale de second terme:
$$E_2^{p,q}=\pi^p(A\otimes \pi^q(\mathcal E_{\infty}(V)))$$
si $\pi^*(A)$ ou $\pi^*(\mathcal E_{\infty}(V)))$ sont finis (finis en chaque dimension) 
celle-ci converge vers 
$\pi^{p+q}(A\amalg_{\tau} \mathcal E_{\infty}(V))$ d\`es que $A$ et $\mathcal E_{\infty}(V)$ sont connexes.
\end{proposition}

\begin{proof}[Preuve]
On filtre l'alg\`ebre $A$ par le degr\'e diff\'erentiel. 
De cette filtration  on d\'eduit une filtration de $A\amalg_{\tau} \mathcal E_{\infty}(V)$.
En effet, la filtration sur $A$ donne une filtration sur $A\otimes B$, et par suite via le morphisme
$p$ une filtration sur $A\amalg\mathcal E_{\infty}(V)$, ce qui permet de construire une suite spectrale.

Afin de d\'eterminer le second terme
de cette suite spectrale, on se ram\`ene au cas d'un coproduit $A\amalg \mathcal E_{\infty}(V)$.
En effet, d'apr\`es l'hypoth\`ese de 1-connexit\'e de $A$, les suites spectrales associ\'ees \`a
$A\amalg_{\tau} \mathcal E_{\infty}(V)$ et $A\amalg \mathcal E_{\infty}(V)$ ont m\^eme second terme.

On remarque que pour un coproduit cette suite spectrale d\'eg\'en\`ere au terme $E_2$. D'apr\`es le
quasi-isomorphisme $A\amalg \mathcal E_{\infty}(V)\longrightarrow A\otimes \mathcal E_{\infty}(V)$,
on en d\'eduit que:
$$E_2^{p,q}=\pi^p(A\otimes \pi^q(\mathcal E_{\infty}(V))).$$
Les alg\`ebres 
$A$ et $\mathcal E_{\infty}(V)$ \'etant suppos\'ees connexes, on obtient une suite spectrale du premier quadrant.
Cette suite spectrale converge vers 
$\pi^{p+q}(A\amalg_{\tau} \mathcal E_{\infty}(V)).$
\end{proof}

{\bf \noindent Le mod\`ele de la fibre}\qua
Soit $p:X\longrightarrow Y$ une fibration de fibre $F$, avec $Y$ suppos\'e $1$-connexe, on suppose que la 
cohomologie de $Y$ ou de $F$ est finie.
On se donne $\mathcal M_Y$ un mod\`ele cofibrant de $\Omega^{\mathcal E_{\infty}}(Y)$. 
On fixe un mod\`ele quasi-libre de
$$\mathcal M_Y \stackrel{\phi_Y}{\longrightarrow}\Omega^{\mathcal E{\infty}}(Y)\stackrel{p^*}{\longrightarrow}
\Omega^{\mathcal E_{\infty}}(X).$$
Soit
$$p':\mathcal M_Y \longrightarrow \mathcal M_Y\amalg_{\tau}\mathcal M_F.$$
On va montrer que la $\mathcal E_{\infty}$-alg\`ebre $\mathcal M_F$ est un mod\`ele cofibrant de $F$.
\begin{diagram}
\Omega^{\mathcal E_{\infty}}(Y)&\rTo^{p^*}&\Omega^{\mathcal E_{\infty}}(X)&\rTo^{i^*}&\Omega^{\mathcal E_{\infty}}(F)\\
\uTo^{\phi_Y}&&\uTo^{\phi_X}&&\uDotsto^{\alpha}&&\\
\mathcal M_Y &\rTo^{p'}& \mathcal M_Y\amalg_{\tau}\mathcal M_F &\rTo^{i'}&\mathcal M_F\\
\end{diagram}
Comme $\mathcal M_F$ est la cofibre de $p'$ il existe un morphisme
$\alpha:\mathcal M_F\longrightarrow \Omega^{\mathcal E_{\infty}}(F)$ compl\'etant le diagramme.

\begin{theorem}
Le morphisme $\alpha:\mathcal M_F\longrightarrow \Omega^{\mathcal E_{\infty}}(F)$
est une fibration triviale.
\end{theorem}
\begin{proof}[Preuve]
On a par d\'efinition
$\alpha i'=i^*\phi_X:\mathcal M_Y\amalg_{\tau}\mathcal M_F\longrightarrow \Omega^{\mathcal E_{\infty}}(F)$.
Or les morphismes $\phi_X$, $i^*$ et $i'$ sont des fibrations, on en d\'eduit que $\alpha$ est aussi une fibration.

Appliquons le foncteur $\Omega^{\mathcal E_{\infty}}$ au diagramme d'ensembles bisimpliciaux:
\begin{diagram}
\beta X &&\rTo^{a}&& S_{.,.}f\\
&\rdTo^{\beta p}&&\ruTo\\
&&\beta{Y}&&\\
\end{diagram}
On obtient:
\begin{diagram}
\Omega^{\mathcal E_{\infty}}(X) &&\lTo^{a^*}&& \Omega^{\mathcal E_{\infty}}_{bi}(S_{.,.}f)\\
&\luTo^{p^*}&&\ldTo\\
&&\Omega^{\mathcal E_{\infty}}(Y)&&\\
\end{diagram}
o\`u le morphime $a^*$ est une fibration acyclique dans la cat\'egorie des $R$-modules diff\'erentiels.
Consid\'erons le diagramme suivant:\eject\vglue-15mm
\begin{diagram}
\Omega^{\mathcal E_{\infty}}(X)&&\rTo^{i^*}&&\Omega^{\mathcal E_{\infty}}(F)\\
\uTo^{\phi_X}&\luTo^{a^*}&&\ruTo&\uTo^{\alpha}\\
&&\Omega^{\mathcal E_{\infty}}_{bi}(S_{.,.}f) &&\\
&\ruDotsto^{\psi}\\
\mathcal M_Y\amalg_{\tau}\mathcal M_F&&\rTo^{i'}&&\mathcal M_F\\
\end{diagram}
L'existence de ${\psi}$ (morphisme dans $\Mdg$) est assur\'ee par le fait que $a^*$ est une fibration
triviale. 

Ce morphisme induit un morphisme de suites spectrales. En particulier, on a:
$$\psi_2^{p,q}:\pi^p(\mathcal M_Y\otimes \pi^q(\mathcal M_F))\longrightarrow H^p(B;H^q(F;R))$$
et
$$\psi_2^{p,0}=\pi(\phi_Y):\pi^p(\mathcal M_Y)\longrightarrow H^p(B;R)$$
$$\psi_2^{0,q}=\pi(\alpha):\pi^q(\mathcal M_F)\longrightarrow H^q(F;R)$$
On sait que $\phi_X$ est un quasi-isomorphisme et que $\psi_2^{p,0}=\pi(\phi_Y)$ est un isomorphisme.
On en d\'eduit que $\psi_2^{0,q}=\pi(\alpha)$ est aussi un isomorphisme par un th\'eor\`eme de
comparaison de suites spectrales (\cite {MC}, section 3.1.1).
\end{proof}

\end{subsection}

\begin{subsection}{Transgression dans la suite spectrale de Leray-Serre}

On fixe $R=k$ un corps. Rappelons la d\'efinition classique de la transgression \cite {EM}.

\begin{definition}
Soit $F \stackrel i \longrightarrow Y \stackrel p \longrightarrow X$ une fibration on consid\`ere le 
diagramme:
\begin{diagram}
H^n(Y;k) & \rTo^{i^*} & H^n(F;k) & \rTo^{\delta} & H^{n+1}(Y,F;k) & \rTo & H^{n+1}(F;k)\\
         &            & \uTo     &               & \uTo^{{p'}^*}     &     
&\uTo^{p^*}\\    &            & H^n(*;k) & \rTo^{\delta} & H^{n+1}(X,*;k)   &
\rTo &H^{n+1}(X;k)\\     \end{diagram}
On dit que $[x]\in H^{n+1}(X,*;k)$ transgresse $[y]\in H^n(F;k)$
si ${p'}^*([x])=\delta ([y])$.
\end{definition}

On suppose que $\pi_1(X)$ op\`ere trivialement sur $H^n(F,k)$.
On rappelle que, pour la suite spectrale de Leray-Serre associ\'ee \`a cette fibration,
on a $E^{0,n}_{n+1}\subset E^{0,n}_2\cong H^n(F,k)$.
Il est \'equivalent de dire que $[x]$ est une
transgression de $[y]\in E^{0,n}_{n+1}$,
si $[x]$ se projette sur $d_{n+1}[y]$ via l'application:
$$H^{n+1}(X,k)=E^{n+1,0}_2\longrightarrow E^{n+1,0}_{n+1}.$$ 
\begin{proposition}
Soit $[x]$ une transgression de $[y]$; pour tout repr\'esentant $y\in Z^n\Omega^{\mathcal E_{\infty}}(F)$ de $[y]$
et tout repr\'esentant $x\in Z^{n+1}\Omega^{\mathcal E_{\infty}}(X,*)$ de $[x]$, il existe
$u\in \Omega^{\mathcal E_{\infty}}(Y)$ tel que $i^*(u)=y$ et $du=q^*(\alpha_0(x))$.
\end{proposition}

\begin{proof}[Preuve]
On consid\`ere le diagramme suivant:
\begin{diagram}
0&\rTo& \Omega^{\mathcal E_{\infty}}(Y,F) &\rTo^{\alpha}&
\Omega^{\mathcal E_{\infty}}(Y)&\rTo^{i^*}&\Omega^{\mathcal E_{\infty}}(F)&\rTo&0\\
 &    &\uTo^{{q'}^*}         &             &\uTo^{q^*}            &         
&\uTo                  &\\ 0&\rTo& \Omega^{\mathcal E_{\infty}}(X,*)
&\rTo^{\alpha_0}&\Omega^{\mathcal E_{\infty}}(X)&\rTo^{i^*}&
\Omega^{\mathcal E_{\infty}}(*)\\
\end{diagram}
Soit $y\in Z^n\Omega^{\mathcal E_{\infty}}(F)$. Comme le morphisme $i^*$ est surjectif, il existe
$u_1\in\Omega^{\mathcal E^n_{\infty}}(Y)$
tel que $i^*(u_1)=y$. Par hypoth\`ese, puisque $[x]$ est une transgression de
$[y]$, on a ${q'}^*([x])=\delta ([y])$. Il existe $\lambda \in Z^{n+1}\mathcal
E_{\infty}(Y,F)$, $\lambda_1\in Z^{n}\Omega^{\mathcal E_{\infty}}(Y,F)$ tels que
${q'}^*(x)=\lambda +d\lambda_1$  et $\alpha(\lambda)=d(u_1)$. On en
d\'eduit $\alpha({q'}^*(x))=\alpha(\lambda)+\alpha(d\lambda_1)$.
C'est-\`a-dire $q^*(\alpha_0(x))=d(u_1+\alpha(\lambda_1))$. Comme on a aussi
$i^*(u_1+\alpha(\lambda_1))=y$, on en conclut que $u=u_1+\alpha(\lambda_1)$ convient.
\end{proof}

{\noindent \bf La transgression alg\'ebrique}\qua
On travaille d\'esormais sur
$\overline {\mathbb F_p}$ la cl\^oture alg\'eb\-rique de $\mathbb F_p$.
On suppose que $A\longrightarrow A\amalg_{\tau} \mathcal O(V)$ est une extension libre, on note
$i:A\longrightarrow A\amalg_{\tau}\mathcal O(V)$ et $p:A\amalg_{\tau}\mathcal O(V)\longrightarrow \mathcal O(V)$
les morphismes canoniques.
On suppose aussi que $A$ est $1$-connexe et que $\mathcal O(V)$ est connexe. On s'int\'eresse \`a la transgression
dans la suite spectrale d\'efinie dans la proposition $4.8$.

\begin{definition}
Un \'el\'ement de $[a]\in \pi^{n+1}(A)$ transgresse $[v]\in E_{n+1}^{0,n}\subset E_{2}^{0,n}\cong\pi^{n}(\mathcal O(V))$,
si $[a]$ se projette sur $d_{n+1}(v)\in E_{n+1}^{n+1,0}$.
\end{definition}

\begin{theorem}[Th\'eor\`eme de Kudo alg\'ebrique]
Si un \'el\'ement $[a]\!\in\!\pi^{n+1}(A)$ est une transgression de $[v]\in \pi^{n}(\mathcal O(V))$, alors
$\mathcal P^s([a])$ transgresse $\mathcal P^s([v])$ et $\beta\mathcal P^s([a])$ transgresse $-\beta \mathcal P^s([v])$.
\end{theorem}

\begin{proof}[Preuve]
La preuve reprend les arguments d\'evelopp\'es par J.P. May dans \cite {May} (th\'eor\`emes $3.3$ et $3.4$).
On applique le th\'eor\`eme $3.3$ de \cite {May} \`a $i$ et $p$. C'est-\`a-dire que l'on construit
un morphisme de ``suspension" $\sigma:Ker(\pi^*(i))\longrightarrow Coker(\pi^*(p))$ en posant
$\sigma([a])=[p(u)]$ tel que $[a]\in Ker(\pi^*(i))$ et $d(u)=i(a)$. Et on montre que ce morphisme
de suspension commute aux op\'erations de Steenrod. On a les formules:
$\sigma \mathcal P^s([a])=\mathcal P^s(\sigma [a])$ et $\sigma \beta \mathcal P^s([a])=-\beta\mathcal P^s(\sigma [a])$.
La deuxi\`eme partie de la preuve est aussi une adaptation ad hoc de la preuve du th\'eor\`eme $3.4$:
si $[v]\in \pi^n(\mathcal O(V))$ est transgressif alors il est repr\'esent\'e par un \'el\'ement $p(u)$ tel que
$d(u)=i(a)$. Le r\'esultat est une cons\'equence des formules de commutation des op\'erations de Steenrod \`a la suspension.
\end{proof}

{\noindent \bf Remarque}\qua le th\'eor\`eme de Kudo classique est alors un corollaire du r\'esultat pr\'ec\'edent appliqu\'e
au mod\`ele de la suite spectrale de Leray-Serre.

\end{subsection}

\end{section}

\renewcommand\refname{R\'ef\'erences}

\Addresses\recd

\end{document}